\documentclass[12pt]{amsart}

\usepackage{amsmath,xspace,amssymb,mathrsfs}
\usepackage{color}

\input xy
\xyoption{all}
\xyoption{2cell}
\UseAllTwocells
\CompileMatrices

\newcommand{\Spec}{\operatorname{Spec}}
\renewcommand{\phi}{\varphi}

\newcommand{\Ker}{\operatorname{Ker}}

\newcommand{\Ima}{\operatorname{Im}}
\newcommand{\Min}{\operatorname{Min}}

\newcommand{\Max}{\operatorname{Max}}

\newcommand{\Sp}{\operatorname{Sp}}
\newcommand{\Ann}{\operatorname{Ann}}

\newtheorem{proposition}{Proposition}[section]
\newtheorem{lemma}[proposition]{Lemma} 
\newtheorem{corollary}[proposition]{Corollary}
\newtheorem{theorem}[proposition]{Theorem}

\newtheorem{prop-def}[proposition]{Proposition and definition}

\theoremstyle{definition}
\newtheorem{definition}[proposition]{Definition}

\newtheorem{remark}[proposition]{Remark}

\begin{document}

\title[Gelfand rings, clean rings and their dual rings]{Characterizations of Gelfand rings specially clean rings and their dual rings}

\author[M. Aghajani and A. Tarizadeh]{Mohsen Aghajani and Abolfazl Tarizadeh}
\address{ Department of Mathematics, Faculty of Basic Sciences,
University of Maragheh \\
P. O. Box 55136-553, Maragheh, Iran.
 }
\email{aghajani14@gmail.com, ebulfez1978@gmail.com}


\begin{abstract}
In this paper, new criteria for zero dimensional rings, Gelfand rings, clean rings and mp-rings are given. A new class of rings is introduced and studied, we call them purified rings. Specially, reduced purified rings are characterized. New characterizations for pure ideals of reduced Gelfand rings and mp-rings are provided. It is also proved that if the topology of a scheme is Hausdorff, then the affine opens of that scheme is stable under taking finite unions. In particular, every compact scheme is an affine scheme. \\
\end{abstract}

\maketitle

\section{Introduction}

This paper is devoted to study two fascinating classes of commutative rings which are so called Gelfand rings and mp-rings. Recall that a ring $A$ is said to be a Gelfand ring (or, pm-ring) if each prime ideal of $A$ is contained in a unique maximal ideal of $A$. Dually, a ring $A$ is called a mp-ring if each prime ideal of $A$ contains a unique minimal prime ideal of $A$. In this paper, we also introduce and study a new class of rings, we call them purified rings. It is shown that purified rings are actually the dual of clean rings. Gelfand rings, and especially clean rings, have been greatly studied in the literature over the past fifty years. Theorem \ref{Theorem VII} improves and generalizes many major results in the literature which are related to the characterization of clean rings. In what follows we shall try to give an account of the new results of this paper in detail. \\
Gelfand rings have been the main subject of many articles in the literature over the years and are still of current interest, see e.g. \cite{Al Ezeh}, \cite{Anderson}, \cite{Contessa}, \cite{Contessa 2}, \cite{Contessa 3}, \cite{Marco-Orsatti} and \cite{McGovern}. The paper \cite{Marco-Orsatti} can be viewed as a starting point of investigations of Gelfand rings in the commutative case. One of the main results of this paper is that ``a ring is a Gelfand ring if and only if its maximal spectrum is the Zariski retraction of its prime spectrum'', see \cite[Theorem 1.2]{Marco-Orsatti}. This result plays a major role in proving some results of our paper. In Theorem \ref{Theorem V}, we prove the dual of this result which states that a ring is a mp-ring if and only if its minimal spectrum is the flat retraction of its prime spectrum. Theorem \ref{Theorem Tariz-Couch} is a further important result which is obtained in this direction. Although mp-rings have been around for some time, there is no substantial account of their characterizations and basic properties in the literature. This may be because Gelfand rings are tied up with the Zariski topology, see Theorem \ref{Theorem III}. By contrast, we show that mp-rings are tied up with the flat topology, see Theorem \ref{Theorem V}. The flat topology is less known than the Zariski topology in the literature. For the flat topology please consider \S 2 and for more details see \cite{Abolfazl}. In this paper, we give a coherent account of mp-rings and their basic and sophisticated properties, see Theorems \ref{Theorem V}, \ref{Theorem Tariz-Couch},
\ref{Theorem VI} and \ref{Theorem 20}. \\
Clean rings, as a subclass of Gelfand rings, are also investigated in this paper. Recall that a ring $A$ is called a clean ring if each element of it can be written as a sum of an idempotent and an invertible element of that ring. This simple definition has some spectacular equivalents, see Theorem \ref{Theorem VII}. This result, in particular, tells us that if $A$ is a clean ring then a system of equations $f_{i}(x_{1},...,x_{n})=0$ with $i=1,...,d$ over $A$ has a solution in $A$ provided that this system has a solution in each local ring $A_{\mathfrak{m}}$ with $\mathfrak{m}$ a maximal ideal of $A$. Clean rings have been extensively studied in the literature over the past and recent years, see e.g. \cite{Anderson}, \cite{Ara et al. 1}, \cite{Couchot}, \cite{Goodearl-Warfield}, \cite{McGovern} and \cite{Nicholson}. But according to \cite{Ara et al. 1}, although examples and constructions of exchange rings abound, there is a pressing need for new constructions to aid the development of the theory (note that in the commutative case, exchange rings and clean rings are the same, see Theorem \ref{Theorem VII}). Toward realizing this purpose, then Theorem \ref{Theorem VII} can be considered as the culmination and strengthening of all of the results in the literature which are related to the characterization of commutative clean rings. \\
Then we introduce and study a new class of rings, we call them purified rings. Specifically, in Theorem \ref{Theorem X} we characterize reduced purified rings. This result, in particular, tells us that if $A$ is a reduced purified ring then a system of equations $f_{i}(x_{1},...,x_{n})=0$ over $A$ has a solution in $A$ provided that this system has a solution in each domain $A/\mathfrak{p}$ with $\mathfrak{p}$ a minimal prime ideal of $A$. \\
In \cite[Theorems 1.8 and 1.9]{Al Ezeh}, the pure ideals of a reduced Gelfand ring are characterized. In Theorem \ref{Theorem 19}, we have improved these results, especially a very simple proof is given for \cite[Theorem 1.8]{Al Ezeh}. Then in Theorem \ref{Theorem 20} and Corollary \ref{Corollary 29}, the pure ideals of a reduced mp-ring are characterized. Corollaries \ref{Corollary 1909} and \ref{Corollary 1910} are further results which are obtained in this direction. \\
Our paper also contains some geometric results. In particular, we show that every compact scheme is an affine scheme, see Theorem \ref{Theorem VIII}. It is also shown that the underlying space of a separated scheme or more generally a quasi-separated scheme is Hausdorff if and only if every point of it is a closed point, see Corollary \ref{Th I}. \\

\section{Preliminaries}

Here we recall some material which is needed in the sequel. In this paper all rings are commutative. An ideal $I$ of a ring $A$ is called a \emph{pure} ideal if the canonical ring map $A\rightarrow A/I$ is a flat ring map.  \\

\begin{theorem}\label{Proposition I} Let $I$ be an ideal of a ring $A$. Then the following statements are equivalent. \\
$\mathbf{(i)}$ $I$ is a pure ideal. \\
$\mathbf{(ii)}$ $\Ann(f)+I=A$ for all $f\in I$. \\
$\mathbf{(iii)}$ If $f\in I$ then there exists some $g\in I$ such that $f(1-g)=0$. \\
$\mathbf{(iv)}$ $I=\{f\in A: \Ann(f)+I=A\}$. \\
\end{theorem}

{\bf Proof.} See \cite[Chap 7, Proposition 2]{Borceux} or \cite[Tag 04PS]{Johan}. $\Box$ \\

If a prime ideal of a ring $A$ is a pure ideal, then it is a minimal prime ideal of $A$. Pure ideals are quite interesting and have important applications even in non-commutative situations, see \cite[Chaps. 7, 8]{Borceux}. \\
If $X$ is a topological space, then by $\pi_{0}(X)$ we mean the set of connected components of $X$ equipped with the quotient topology. A subspace $Y$ of a topological space $X$ is called a \emph{retraction} of $X$ if there exists a continuous map $\gamma:X\rightarrow Y$ such that $\gamma(y)=y$ for all $y\in Y$. Such a map $\gamma$ is called a retraction map. A quasi-compact and Hausdorff topological space is called a compact space. \\
A subset $E\subseteq\Spec(A)$ is said to be stable under the generalization if whenever $\mathfrak{p}\in E$ and $\mathfrak{q}$ is a prime ideal of $A$ such that $\mathfrak{q}\subseteq\mathfrak{p}$, then $\mathfrak{q}\in E$. The dual notion is called stable under the specialization. If $I$ is a pure ideal of a ring $A$ then by Theorem \ref{Proposition I}, $V(I)$ is stable under the generalization. \\
If an ideal of a ring $A$ is generated by a set of idempotents of $A$, then it is called a \emph{regular ideal} of $A$. Every regular ideal is a pure ideal, but the converse does not necessarily hold. Every maximal element of the set of proper and regular ideals of $A$ is called a \emph{max-regular} ideal of $A$. The set of max-regular ideals of $A$ is called the \emph{Pierce spectrum} of $A$ and denoted by $\Sp(A)$. It is a compact and totally disconnected space whose basis opens are precisely of the form $U_{f}=\{M\in\Sp(A): f\notin M\}$ where $f\in A$ is an idempotent. We also have the homeomorphism $\pi_{0}\big(\Spec(A)\big)\simeq\Sp(A)$, for more details see \cite{Abolfazl}. \\
If $A$ is a ring, then there exists a (unique) topology over $\Spec(A)$ such that the collection of subsets $V(f)=\{\mathfrak{p}\in\Spec(A): f\in\mathfrak{p}\}$ with $f\in A$ forms a sub-basis for the opens of this topology. It is called the flat (or, inverse) topology. Moreover there is a (unique) topology over $\Spec(A)$ such that the collection of subsets $D(f)\cap V(g)$ with $f,g\in A$ forms a sub-basis for the opens of this topology. It is called the patch (or, constructible) topology. If $\mathfrak{p}$ is a prime ideal of $A$, then $\Lambda(\mathfrak{p})=\{\mathfrak{q}\in\Spec(A): \mathfrak{q}\subseteq\mathfrak{p}\}$ is the closure of $\{\mathfrak{p}\}$ in $\Spec(A)$ with respect to the flat topology. It can be shown that the flat closed subsets of $\Spec(A)$ are precisely of the form $\Ima\phi^{\ast}$ where $\phi:A\rightarrow B$ is a flat ring map.
For more details on the flat and patch topologies please consider \cite{Abolfazl}. \\
The following result is due to Grothendieck and has found interesting applications in this paper. \\

\begin{theorem}\label{Theorem II1} The map $f\rightsquigarrow D(f)$ is a bijective function from the set of idempotents of a ring $A$ onto the set of clopen $($both open and closed$)$ subsets of $\Spec(A)$. \\
\end{theorem}

{\bf Proof.} See \cite[00EE]{Johan}. $\Box$ \\

Let $\phi:A\rightarrow B$ be a morphism of rings. We say that the idempotents of $A$ can be lifted along $\phi$ if $e\in B$ is an idempotent, then there exists an idempotent $e'\in A$ such that $\phi(e')=e$. \\

\begin{theorem}\label{Theorem IX} Let $X$ be a compact and totally disconnected topological space. Then the set of clopens of $X$ forms a basis for the topology of $X$. If moreover, $X$ has an open covering $\mathscr{C}$ with the property that every open subset of each member of $\mathscr{C}$ is again a member of $\mathscr{C}$, then there exist a finite number $W_{1},...,W_{q}\in\mathscr{C}$ of pairwise disjoint clopens of $X$ such that $X=\bigcup\limits_{k=1}^{q}W_{k}$. \\
\end{theorem}

{\bf Proof.} The first part of the assertion is deduced from \cite[Theorem 6.1.23]{Engelking}, and the second part follows from the first part. $\Box$ \\

\section{Maximality of primes}

\begin{lemma}\label{Lemma IV} Let $S$ and $T$ be multiplicative subsets of a ring $A$. Then $S^{-1}A\otimes_{A}T^{-1}A=0$ iff there exist $f\in S$ and $g\in T$ such that $fg=0$. \\
\end{lemma}

{\bf Proof.} Suppose $S^{-1}A\otimes_{A}T^{-1}A=0$. Setting $M:=S^{-1}A$. Then $T^{-1}M\simeq M\otimes_{A}T^{-1}A=0$. Thus the image of the unit of $M$ under the canonical map $M\rightarrow T^{-1}M$ is zero. Hence, there exists some $g\in T$ such that $g/1=0$ in $M$. Thus there is some $f\in S$ such that $fg=0$. Conversely, if $r/s\otimes r'/s'$ is a pure tensor of $S^{-1}A\otimes_{A}T^{-1}A$ then we may write $r/s\otimes r'/s'=rf/sf\otimes r'g/s'g=r/sf\otimes r'fg/s'g=0$.
$\Box$ \\

\begin{lemma}\label{Lemma III} If $\mathfrak{p}$ and $\mathfrak{q}$ are distinct minimal prime ideals of a ring $A$, then there exist $f\in A\setminus\mathfrak{p}$ and $g\in A\setminus\mathfrak{q}$ such that $fg=0$. \\
\end{lemma}

{\bf Proof.} It suffices to show that $0\in S:=(A\setminus\mathfrak{p})(A\setminus\mathfrak{q})$. If $0\notin S$ then there exists a prime ideal $\mathfrak{p}'$ of $A$ such that $\mathfrak{p}'\cap S=\emptyset$. This yields that $\mathfrak{p}=\mathfrak{p}'=\mathfrak{q}$. But this is a contradiction. $\Box$ \\

Let $A$ be a ring. Consider the relation $S=\{(\mathfrak{p},\mathfrak{q})\in X^{2}: A_{\mathfrak{p}}\otimes_{A}A_{\mathfrak{q}}\neq0\}$ on $X=\Spec(A)$. By Lemma \ref{Lemma IV}, $(\mathfrak{p},\mathfrak{q})\in S$ iff $0\notin(A\setminus\mathfrak{p})(A\setminus\mathfrak{q})$. This relation is reflexive and symmetric.
Let $\sim_{S}$ be the equivalence relation generated by $S$. Thus $\mathfrak{p}\sim_{S}\mathfrak{q}$ if and only if there exists a finite set $\{\mathfrak{p}_{1},...,\mathfrak{p}_{n}\}$ of prime ideals of $A$ with $n\geqslant2$ such that $\mathfrak{p}_{1}=\mathfrak{p}$, $\mathfrak{p}_{n}=\mathfrak{q}$ and $A_{\mathfrak{p}_{i}}\otimes_{A}A_{\mathfrak{p}_{i+1}}\neq0$ for all $1\leqslant i\leqslant n-1$. Note that it may happen that $\mathfrak{p}\sim_{S}\mathfrak{q}$ but $ A_{\mathfrak{p}}\otimes_{A}A_{\mathfrak{q}}=0$. Also note that in Theorem \ref{Theorem II} (x), by $[\mathfrak{p}]$ we mean the equivalence class of $\sim_{S}$ containing $\mathfrak{p}$. \\
In the following result, the criteria (i)-(iv) are well known, and the remaining are new. The Zariski, flat and patch topologies on $\Spec(A)$ are denoted by $\mathcal{Z}$, $\mathcal{F}$ and $\mathcal{P}$, respectively. \\

\begin{theorem}\label{Theorem II} For a ring $A$ the following statements are equivalent. \\
$\mathbf{(i)}$ Every prime ideal of $A$ is a maximal ideal. \\
$\mathbf{(ii)}$ $\mathcal{Z}$ is Hausdorff.\\
$\mathbf{(iii)}$ $\mathcal{Z}=\mathcal{P}$. \\
$\mathbf{(iv)}$ $A/\mathfrak{N}$ is absolutely flat where $\mathfrak{N}$ is the nil-radical of $A$. \\
$\mathbf{(v)}$ If $\mathfrak{p}$ and $\mathfrak{q}$ are distinct prime ideals of $A$, then there exist $f\in A\setminus\mathfrak{p}$ and $g\in A\setminus\mathfrak{q}$ such that $fg=0$. \\
$\mathbf{(vi)}$ $\mathcal{F}$ is Hausdorff. \\
$\mathbf{(vii)}$ Every flat epimorphism of rings with source $A$ is surjective. \\
$\mathbf{(viii)}$ If $\mathfrak{p}$ is a prime ideal of $A$, then the canonical map $\pi_{\mathfrak{p}}:A\rightarrow A_{\mathfrak{p}}$ is surjective. \\
$\mathbf{(ix)}$ $\mathcal{Z}=\mathcal{F}$. \\
$\mathbf{(x)}$ If $\mathfrak{p}$ is a prime ideal of $A$, then $[\mathfrak{p}]=\{\mathfrak{p}\}$. \\
\end{theorem}

{\bf Proof.} For $\mathbf{(i)}\Leftrightarrow\mathbf{(ii)}$
see \cite[Theorem 3.6]{Gilmer}, for $\mathbf{(i)}\Leftrightarrow\mathbf{(iii)}$ see \cite[Theorem 2.1]{Chapman et al}, and for $\mathbf{(i)}\Leftrightarrow\mathbf{(iv)}$ see \cite[Theorem 3.1]{Huckaba}. \\
$\mathbf{(i)}\Rightarrow\mathbf{(v)}:$ It follows from Lemma \ref{Lemma III}.  \\
$\mathbf{(v)}\Rightarrow\mathbf{(i)}:$ It is straightforward. \\
$\mathbf{(i)}\Rightarrow\mathbf{(vi)}:$ If $\mathfrak{p}$ and $\mathfrak{q}$ are distinct prime ideals of $A$ then by the hypothesis, $\mathfrak{p}+\mathfrak{q}=A$. Thus there are $f\in\mathfrak{p}$ and $g\in\mathfrak{q}$ such that $f+g=1$. Therefore $V(f)\cap V(g)=\emptyset$. \\
$\mathbf{(vi)}\Rightarrow\mathbf{(i)}:$ If $\mathfrak{m}$ is a maximal ideal of $A$, then $\Lambda(\mathfrak{m})=\{\mathfrak{m}\}$, because in a Hausdorff space each point is a closed point. This shows that every prime ideal of $A$ is a maximal ideal. \\
$\mathbf{(i)}\Rightarrow\mathbf{(vii)}:$ Let $\phi:A\rightarrow B$ be a flat epimorphism of rings. It suffices to show that the induced morphism $\phi_{\mathfrak{p}}:A_{\mathfrak{p}}\rightarrow B_{\mathfrak{p}}$ is surjective for all $\mathfrak{p}\in\Spec(A)$. If $\mathfrak{p}B\neq B$, then we show that $\phi_{\mathfrak{p}}$ is an isomorphism of rings. Clearly $\phi_{\mathfrak{p}}$ is a flat epimorphism of rings, because both flat ring maps and epimorphisms of rings are stable under base change. Also, $B_{\mathfrak{p}}$ as $A_{\mathfrak{p}}-$module is faithfully flat, since $\mathfrak{p}B\neq B$. Hence, $\phi_{\mathfrak{p}}$ is an isomorphism of rings, since it is well known that every faithfully flat epimorphism of rings is an isomorphism. But if $\mathfrak{p}B=B$, then we show that $B_{\mathfrak{p}}\simeq A_{\mathfrak{p}}\otimes_{A}B=0$. If $A_{\mathfrak{p}}\otimes_{A}B\neq0$ then it has a prime ideal $P$, and so in the following pushout diagram: $$\xymatrix{
A\ar[r]^{\phi}\ar[d]^{\pi}&B
\ar[d]^{\mu}\\A_{\mathfrak{p}}\ar[r]^{\lambda\:\:\:\:\:\:\:\:\:}&
A_{\mathfrak{p}}\otimes_{A}B}$$ we have $\lambda^{-1}(P)=\mathfrak{p}A_{\mathfrak{p}}$, since by the hypothesis every prime ideal of $A$ is a maximal ideal. Thus $\mathfrak{p}=\phi^{-1}(\mathfrak{q})$ where $\mathfrak{q}:=\mu^{-1}(P)$. It follows that $\mathfrak{p}B\subseteq\mathfrak{q}\neq B$ which is a contradiction. \\
$\mathbf{(vii)}\Rightarrow\mathbf{(viii)}:$ There is nothing to prove. \\
$\mathbf{(viii)}\Rightarrow\mathbf{(i)}:$ For each $f\in A\setminus\mathfrak{p}$ there exists some $g\in A$ such that $g/1=1/f$ in $A_{\mathfrak{p}}$. Thus there exists an element $h\in A\setminus\mathfrak{p}$ such that $h(1-fg)=0$. It follows that $1-fg\in\mathfrak{p}$ and so $A/\mathfrak{p}$ is a field. \\
$\mathbf{(vi)}\Rightarrow\mathbf{(ix)}:$ If $X=\Spec(A)$ then the map $(X,\mathcal{P})\rightarrow(X,\mathcal{F})$ given by $x\rightsquigarrow x$ is continuous, since the patch topology is finer than the flat topology. It is also a closed map, because the space $(X,\mathcal{P})$ is quasi-compact and $(X,\mathcal{F})$ is Hausdorff. Hence, it is a homeomorphism. Thus $\mathcal{F}=\mathcal{P}$. Using (iii), then we get that $\mathcal{Z}=\mathcal{F}$. \\
$\mathbf{(ix)}\Rightarrow\mathbf{(i)}:$ If $\mathfrak{p}$ is a prime of $A$, then $V(\mathfrak{p})=\Lambda(\mathfrak{p})$ and so $\mathfrak{p}$ is a maximal ideal. \\
$\mathbf{(i)}\Rightarrow\mathbf{(x)}:$ Let $\mathfrak{m}$ be a maximal ideal of $A$ and $\mathfrak{m}'\in[\mathfrak{m}]$. Thus there exists a finite set $\{\mathfrak{m}_{1},...,\mathfrak{m}_{n}\}$ of maximal ideals of $A$ with $n\geqslant2$ such that $\mathfrak{m}_{1}=\mathfrak{m}$, $\mathfrak{m}_{n}=\mathfrak{m}'$ and $A_{\mathfrak{m}_{i}}\otimes_{A}A_{\mathfrak{m}_{i+1}}\neq0$ for all $1\leqslant i\leqslant n-1$. Thus by Lemma \ref{Lemma III}, $\mathfrak{m}=\mathfrak{m}_{1}=...=\mathfrak{m}_{n}=\mathfrak{m}'$. \\
$\mathbf{(x)}\Rightarrow\mathbf{(i)}:$ Let $\mathfrak{p}$ be a prime of $A$. There is a maximal ideal $\mathfrak{m}$ of $A$ such that $\mathfrak{p}\subseteq\mathfrak{m}$. By Lemma \ref{Lemma IV}, $A_{\mathfrak{p}}\otimes_{A}A_{\mathfrak{m}}\neq0$. Thus $\mathfrak{m}\in[\mathfrak{p}]$ and so $\mathfrak{p}=\mathfrak{m}$.
$\Box$ \\

\begin{theorem}\label{Theorem VIII} If the topology of a scheme $X$ is Hausdorff, then every finite union of affine opens of $X$ is an affine open. In particular, every compact scheme is an affine scheme. \\
\end{theorem}

{\bf Proof.} Using the induction, then it suffices to prove the assertion for two cases. Hence, let $U$ and $V$ be two affine opens of $X$. Every affine open of $X$ is closed, because in a Hausdorff space each quasi-compact subset is closed. It follows that $W:=U\cap V$ is a clopen of $U$. Therefore $W$, $U\setminus W$ and $V\setminus W$ are affine opens, see Theorem \ref{Theorem II1}. It is well known that if a scheme can be written as the disjoint union of a finite number of affine opens, then it is an affine scheme. Hence, $U\cup V$ is an affine open. $\Box$ \\

We use the above theorems to obtain more geometric results. \\

\begin{corollary}\label{Corollary IV} The category of compact $($affine$)$ schemes is anti-equivalent to the category of zero dimensional rings. $\Box$ \\
\end{corollary}

\begin{corollary}\label{Th I} Let $X$ be a scheme which has an affine open covering such that the intersection of any two elements of this covering is quasi-compact. Then the underlying space of $X$ is Hausdorff if and only if every point of $X$ is a closed point. \\
\end{corollary}

{\bf Proof.} The implication ``$\Rightarrow$" is obvious since each point of a Hausdorff space is a closed point. Conversely, if $X=\Spec(A)$ is an affine scheme then every prime of $A$ is a maximal ideal. Thus by Theorem \ref{Theorem II} (ii), $\Spec(A)$ is Hausdorff. For the general case, let $x$ and $y$ be two distinct points of $X$. By the hypothesis, there exist affine opens $U$ and $V$ of $X$ such that $x\in U$, $y\in V$ and $U\cap V$ is quasi-compact. If either $x\in V$ or $y\in U$ then the assertion holds. Because, by what we have proved above, every affine open of $X$ is Hausdorff. Therefore we may assume that $x\notin V$ and $y\notin U$. But $W:=V\setminus(U\cap V)$ is an open subset of $X$ because every quasi-compact (=compact) subset of a Hausdorff space is closed. Clearly $y\in W$ and $U\cap W=\emptyset$. $\Box$ \\

\begin{remark} The hypothesis of Corollary \ref{Th I} is not limitative at all. Because a separated scheme or more generally a quasi-separated scheme has this property, see \cite[Proposition 3.6]{Liu} or \cite[Ex. 4.3]{Hartshorne} for the separated case and \cite[Tag 054D]{Johan} for the quasi-separated case.\\
\end{remark}

\section{Gelfand rings}

\begin{lemma}\label{Lemma V} Let $S$ be a multiplicative subset of a ring $A$. Then the canonical morphism $\pi:A\rightarrow S^{-1}A$ is surjective if and only if $\Ima\pi^{\ast}=\{\mathfrak{p}\in\Spec(A):\mathfrak{p}\cap S=\emptyset\}$ is a Zariski closed subset of $\Spec(A)$. \\
\end{lemma}

{\bf Proof.} The map $\pi^{\ast}:\Spec(S^{-1}A)\rightarrow\Spec(A)$ is a homeomorphism onto its image. If $\Ima\pi^{\ast}$ is Zariski closed, then $\pi^{\ast}$ is a closed map. If $\mathfrak{q}\in\Spec(S^{-1}A)$, then the morphism $A_{\mathfrak{p}}\rightarrow(S^{-1}A)_{\mathfrak{q}}$ induced by $\pi$ is an isomorphism where $\mathfrak{p}=\pi^{-1}(\mathfrak{q})$. Therefore the morphism $(\pi^{\ast},\pi^{\sharp}):\Spec(S^{-1}A)\rightarrow\Spec(A)$ is a closed immersion of schemes. It is well known that a morphism of rings $A\rightarrow B$ is surjective if and only if the induced morphism $\Spec(B)\rightarrow\Spec(A)$ is a closed immersion of schemes. Hence, $\pi$ is surjective. Conversely, if $\pi$ is surjective then $\Ima\pi^{\ast}=V(\Ker\pi)$. $\Box$ \\

If $\mathfrak{p}$ is a prime ideal of a ring $A$, then the image of each $f\in A$ under the canonical ring map $\pi_{\mathfrak{p}}:A\rightarrow A_{\mathfrak{p}}$ is also denoted by $f_{\mathfrak{p}}$. \\

\begin{remark}\label{Remark VIII} Let $A$ be a ring and consider the following relation $R=\{(\mathfrak{p},\mathfrak{q})\in X^{2}:\mathfrak{p}+\mathfrak{q}\neq A\}$ on $X=\Spec(A)$. Clearly it is reflexive and symmetric. Let $\sim_{R}$ be the equivalence relation generated by $R$. Then  $\mathfrak{p}\sim_{R}\mathfrak{q}$ if and only if there exists a finite set $\{\mathfrak{p}_{1},...,\mathfrak{p}_{n}\}$ of prime ideals of $A$ with $n\geqslant2$ such that $\mathfrak{p}_{1}=\mathfrak{p}$, $\mathfrak{p}_{n}=\mathfrak{q}$ and $\mathfrak{p}_{i}+\mathfrak{p}_{i+1}\neq A$ for all $1\leqslant i\leqslant n-1$. Note that it may happen that $\mathfrak{p}\sim_{R}\mathfrak{q}$ but $\mathfrak{p}+\mathfrak{q}=A$. Obviously  $\Spec(A)/\sim_{R}=\{[\mathfrak{m}]: \mathfrak{m}\in\Max A\}=\{[\mathfrak{p}]:\mathfrak{p}\in\Min A\}$. It is important to notice that in Theorem \ref{Theorem III} (xi) by $\Spec(A)/\sim_{R}$ we mean the quotient of the ``Zariski'' space $\Spec(A)$ modulo $\sim_{R}$. By contrast, in Theorem \ref{Theorem V} (viii), by $\Spec(A)/\sim_{R}$ we mean the quotient of the ``flat'' space $\Spec(A)$ modulo $\sim_{R}$. Also note that in Theorems \ref{Theorem III} and \ref{Theorem V}, by $[\mathfrak{p}]$ we mean the equivalence class of $\sim_{R}$ containing $\mathfrak{p}$. \\
\end{remark}

Recall that a ring $A$ is called a Gelfand ring if each prime ideal of $A$ is contained in a unique maximal ideal of $A$. In the following result, the criteria (i)-(v) are well known, and the remaining are new. \\

\begin{theorem}\label{Theorem III} For a ring $A$ the following statements are equivalent. \\
$\mathbf{(i)}$ $A$ is a Gelfand ring. \\
$\mathbf{(ii)}$ $\Max(A)$ is the Zariski retraction of $\Spec(A)$. \\
$\mathbf{(iii)}$ $\Spec(A)$ is a normal space with respect to the Zariski topology.  \\
$\mathbf{(iv)}$ If $f\in A$, then there exist $g,h\in A$ such that $(1+fg)(1+f'h)=0$ where $f'=1-f$. \\
$\mathbf{(v)}$ If $\mathfrak{m}$ and $\mathfrak{m}'$ are distinct maximal ideals of $A$, then there exist $f\in A\setminus\mathfrak{m}$ and $g\in A\setminus\mathfrak{m}'$ such that $fg=0$. \\
$\mathbf{(vi)}$ If $\mathfrak{m}$ is a maximal ideal of $A$, then the canonical map $\pi_{\mathfrak{m}}:A\rightarrow A_{\mathfrak{m}}$ is surjective. \\
$\mathbf{(vii)}$ If $\mathfrak{m}$ is a maximal ideal of $A$, then $[\mathfrak{m}]=\Lambda(\mathfrak{m})$. \\
$\mathbf{(viii)}$ If $\mathfrak{m}$ is a maximal ideal of $A$, then $\Lambda(\mathfrak{m})$ is a Zariski closed subset of $\Spec(A)$. \\
$\mathbf{(ix)}$ If $\mathfrak{m}$ and $\mathfrak{m}'$ are distinct maximal ideals of $A$, then $\Ker\pi_{\mathfrak{m}}+\Ker\pi_{\mathfrak{m}'}=A$. \\
$\mathbf{(x)}$ If $\mathfrak{m}$ and $\mathfrak{m}'$ are distinct maximal ideals of $A$, then there exists some $f\in A$ such that $f_{\mathfrak{m}}=0$ and $f_{\mathfrak{m}'}=1$. \\
$\mathbf{(xi)}$ The map $\eta:\Max(A)\rightarrow\Spec(A)/\sim_{R}$ given by $\mathfrak{m}\rightsquigarrow[\mathfrak{m}]$ is a homeomorphism. \\
$\mathbf{(xii)}$ If $\mathfrak{m}$ is a maximal ideal of $A$, then $\Lambda(\mathfrak{m})=V(\Ker\pi_{\mathfrak{m}})$. \\
\end{theorem}

{\bf Proof.} For $\mathbf{(i)}\Leftrightarrow\mathbf{(ii)}\Leftrightarrow\mathbf{(iii)}$
see \cite[Theorem 1.2]{Marco-Orsatti}, and for $\mathbf{(i)}\Leftrightarrow\mathbf{(iv)}\Leftrightarrow\mathbf{(v)}$ see \cite[Theorem 4.1]{Contessa} and \cite[Proposition 1.1]{Contessa 2}. \\
$\mathbf{(v)}\Rightarrow\mathbf{(vi)}:$ It suffices to show that the morphism $\phi_{\mathfrak{m}'}:A_{\mathfrak{m}'}\rightarrow
(A_{\mathfrak{m}})_{\mathfrak{m}'}$ induced by $\phi:=\pi_{\mathfrak{m}}$ is surjective for all $\mathfrak{m}'\in\Max(A)$. Clearly $\phi_{\mathfrak{m}}$ is an isomorphism. If $\mathfrak{m}'\neq\mathfrak{m}$, then by Lemma \ref{Lemma IV}, $(A_{\mathfrak{m}})_{\mathfrak{m}'}\simeq A_{\mathfrak{m}}\otimes_{A}A_{\mathfrak{m}'}=0$. \\
$\mathbf{(vi)}\Rightarrow\mathbf{(v)}:$ Choose some $h\in\mathfrak{m}'\setminus\mathfrak{m}$ then there exists some $a\in A$ such that $1/h=a/1$ in $A_{\mathfrak{m}}$. Thus there exists some $f\in A\setminus\mathfrak{m}$ such that $f(ah-1)=0$. Clearly $g:=ah-1\in A\setminus\mathfrak{m}'$. \\
$\mathbf{(i)}\Rightarrow\mathbf{(vii)}:$ Let $\mathfrak{p}\in[\mathfrak{m}]$. There exists a maximal ideal $\mathfrak{m}'$ of $A$ such that $\mathfrak{p}\subseteq\mathfrak{m}'$. It follows that $\mathfrak{m}\sim_{R}\mathfrak{m}'$. Thus there exists a finite set $\{\mathfrak{p}_{1},...,\mathfrak{p}_{n}\}$ of prime ideals of $A$ with $n\geqslant2$ such that $\mathfrak{p}_{1}=\mathfrak{m}$, $\mathfrak{p}_{n}=\mathfrak{m}'$ and $\mathfrak{p}_{i}+\mathfrak{p}_{i+1}\neq A$ for all $1\leqslant i\leqslant n-1$. By induction on $n$ we shall prove that $\mathfrak{m}=\mathfrak{m}'$. If $n=2$ then $\mathfrak{m}+\mathfrak{m}'\neq A$ and so $\mathfrak{m}=\mathfrak{m}'$. Assume that $n>2$. We have $\mathfrak{p}_{n-2}+\mathfrak{p}_{n-1}\neq A$ and $\mathfrak{p}_{n-1}\subseteq\mathfrak{m}'$. Thus by the hypothesis, $\mathfrak{p}_{n-2}\subseteq\mathfrak{m}'$ and so $\mathfrak{p}_{n-2}+\mathfrak{m}'\neq A$. Thus in the equivalence $\mathfrak{m}\sim_{R}\mathfrak{m}'$ the number of the involved primes is reduced to $n-1$. Therefore by the induction hypothesis, $\mathfrak{m}=\mathfrak{m}'$. \\
$\mathbf{(vii)}\Rightarrow\mathbf{(i)}:$ Let $\mathfrak{p}$ be a prime of $A$ such that $\mathfrak{p}\subseteq\mathfrak{m}$ and  $\mathfrak{p}\subseteq\mathfrak{m}'$ for some maximal ideals $\mathfrak{m}$ and $\mathfrak{m}'$ of $A$. It follows that $\mathfrak{m}\sim_{R}\mathfrak{m}'$ and so $[\mathfrak{m}]=[\mathfrak{m}']$. Thus by the hypothesis, $\mathfrak{m}=\mathfrak{m}'$. \\
$\mathbf{(vi)}\Leftrightarrow\mathbf{(viii)}:$ See Lemma \ref{Lemma V}. \\
$\mathbf{(vi)}\Rightarrow\mathbf{(ix)}:$ We have $V(\Ker\pi_{\mathfrak{m}})=\{\mathfrak{p}\in\Spec(A): \mathfrak{p}\subseteq\mathfrak{m}\}$. If  $\Ker\pi_{\mathfrak{m}}+\Ker\pi_{\mathfrak{m}'}\neq A$, then there exists a maximal ideal $\mathfrak{m}''$ such that $\Ker\pi_{\mathfrak{m}}+\Ker\pi_{\mathfrak{m}'}
\subseteq\mathfrak{m}''$.  This yields that $\mathfrak{m}=\mathfrak{m}''=\mathfrak{m}'$ which is a contradiction.\\
$\mathbf{(ix)}\Rightarrow\mathbf{(v)}:$ There are $f'\in\Ker\pi_{\mathfrak{m}}$ and $g'\in\Ker\pi_{\mathfrak{m}'}$ such that $f'+g'=1$. Thus there exist $f\in A\setminus\mathfrak{m}$ and $g\in A\setminus\mathfrak{m}'$ such that $ff'=gg'=0$. It follows that $fg=0$. \\
$\mathbf{(ix)}\Rightarrow\mathbf{(x)}:$ There exists some $f\in\Ker\pi_{\mathfrak{m}}$ such that $1-f\in\Ker\pi_{\mathfrak{m}'}$. \\
$\mathbf{(x)}\Rightarrow\mathbf{(v)}:$ There exist $g\in A\setminus\mathfrak{m}$ and $h\in A\setminus\mathfrak{m}'$ such that $fg=0$ and $(1-f)h=0$. It follows that $gh=0$. \\
$\mathbf{(i)}\Rightarrow\mathbf{(xi)}:$ For any ring $A$ the map $\eta$ is continuous and surjective, see Remark \ref{Remark VIII}. By (vii), it is injective. It remains to show that its converse $\mu$ is continuous. By (ii), $\gamma=\mu\circ\pi$ is continuous where $\pi:\Spec(A)\rightarrow\Spec(A)/\sim_{R}$ is the canonical morphism. Therefore $\mu$ is continuous. \\
$\mathbf{(xi)}\Rightarrow\mathbf{(i)}:$ It is proved exactly like the implication (vii)$\Rightarrow$(i). \\
$\mathbf{(ix)}\Rightarrow\mathbf{(xii)}:$ Clearly $\Lambda(\mathfrak{m})\subseteq V(\Ker\pi_{\mathfrak{m}})$. Let $\mathfrak{p}$ be a prime ideal of $A$ such that $\Ker\pi_{\mathfrak{m}}\subseteq\mathfrak{p}$. If $\mathfrak{p}\nsubseteq\mathfrak{m}$, then there exists a maximal ideal $\mathfrak{m}'$ of $A$ such that $\mathfrak{p}\subseteq\mathfrak{m}'$. Using the hypothesis, then we get that $\mathfrak{m}'=A$ which is a contradiction.
$\mathbf{(xii)}\Rightarrow\mathbf{(viii)}:$ There is nothing to prove. $\Box$ \\

\begin{remark} Clearly Gelfand rings are stable under taking quotients, and mp-rings are stable under taking localizations. But the localization of a Gelfand ring is not necessarily a Gelfand ring. For example, consider the prime ideals $\mathfrak{p}=(x/1)$ and $\mathfrak{q}=(y/1)$ in $A=(k[x,y])_{P}$ where $k$ is a domain and $P=(x,y)$. Then $A$ is a Gelfand ring but $S^{-1}A$ is not a Gelfand ring where $S=A\setminus\mathfrak{p}\cup\mathfrak{q}$. Dually, $k[x,y]$ is a mp-ring but the quotient $A=k[x,y]/I$ is not a mp-ring, because $\mathfrak{p}=(x+I)$ and $\mathfrak{q}=(y+I)$ are two distinct minimal primes of $A$ which are contained in the prime ideal $(x+I,y+I)$ where $I=(xy)$. \\
\end{remark}

\begin{remark}\label{Remark III} Note that the main result of \cite[Theorem 2.11 (iv)]{Simmons} by Harold Simmons
is not true and the gap is not repairable. It claims that if $\Max(A)$ is Zariski Hausdorff, then $A$ is a Gelfand ring. As a counterexample for the claim, let $p$ and $q$ be two distinct prime numbers, $S=\mathbb{Z}\setminus(p\mathbb{Z}\cup q\mathbb{Z})$ and $A=S^{-1}\mathbb{Z}$. Then $\Max(A)=\{S^{-1}(p\mathbb{Z}),S^{-1}(q\mathbb{Z})\}$ is Zariski Hausdorff because $D(p/1)\cap\Max A=\{S^{-1}(q\mathbb{Z})\}$ and $D(q/1)\cap\Max(A)=\{S^{-1}(p\mathbb{Z})\}$. But $A$ is not a Gelfand ring, since it is nonlocal domain. See also \cite{Simmons 2}, Simmons' erratum for \cite{Simmons}. In Proposition \ref{Proposition II}, we correct his mistake. \\
\end{remark}

\begin{proposition}\label{Proposition II} If $\Max(A)$ is Zariski Hausdorff and $\mathfrak{N}=\mathfrak{J}$, then $A$ is a Gelfand ring where $\mathfrak{J}$ is the Jacobson radical of $A$. \\
\end{proposition}

{\bf Proof.} Let $\mathfrak{m}$ and $\mathfrak{m}'$ be distinct maximal ideals of $A$ both containing a prime $\mathfrak{p}$ of $A$. By the hypotheses, there are $f\in A\setminus\mathfrak{m}$ and $g\in A\setminus\mathfrak{m}'$ such that $\big(D(f)\cap\Max(A)\big)\cap\big(D(g)\cap\Max(A)\big)=\emptyset$. It follows that $fg\in\mathfrak{J}$. Thus there exists a natural number $n\geqslant1$ such that $f^{n}g^{n}=0$. Then either $f\in\mathfrak{p}$ or $g\in\mathfrak{p}$. This is a contradiction, hence $A$ is a Gelfand ring. $\Box$ \\

\begin{corollary}\label{Corollary II} Let $A$ be a ring. Then $\Max(A)$ is Zariski Hausdorff if and only if $A/\mathfrak{J}$ is a Gelfand ring. $\Box$ \\
\end{corollary}

\section{Clean rings}

In this section we give new characterizations for clean rings. \\
By a system of equations over a ring $A$ we mean a finite number of equations $f_{i}(x_{1},...,x_{n})=0$ with $i=1,...,d$ where each $f_{i}(x_{1},...,x_{n})\in A[x_{1},...,x_{n}]$. We say that this system has a solution in $A$ if there exists an n-tuple $(c_{1},...,c_{n})\in A^{n}$ such that $f_{i}(c_{1},...,c_{n})=0$ for all $i$. \\
Remember that a ring $A$ is called a clean ring if each $f\in A$
can be written as $f=e+u$ where $e\in A$ is an idempotent and $u\in A$ is invertible in $A$. In the following result, the criteria (i)-(vii) are well known, and the remaining are new. \\

\begin{theorem}\label{Theorem VII} For a ring $A$ the following statements are equivalent. \\
$\mathbf{(i)}$ $A$ is a clean ring. \\
$\mathbf{(ii)}$ $A$ is an exchange ring, i.e., for each $f\in A$ there exists an idempotent $e\in A$ such that $e\in Af$ and $1-e\in A(1-f)$. \\
$\mathbf{(iii)}$ $A$ is a Gelfand ring and every pure ideal of $A$ is a regular ideal. \\
$\mathbf{(iv)}$ The collection of $D(e)\cap\Max(A)$ with $e\in A$ an idempotent forms a basis for the induced Zariski topology on $\Max(A)$. \\
$\mathbf{(v)}$ The idempotents of $A$ can be lifted modulo each ideal of $A$. \\
$\mathbf{(vi)}$ $A$ is a Gelfand ring and $\Max(A)$ is totally disconnected with respect to the Zariski topology. \\
$\mathbf{(vii)}$ $A$ is a Gelfand ring and for each maximal ideal $\mathfrak{m}$ of $A$, $\Ker\pi_{\mathfrak{m}}$ is a regular ideal of $A$. \\
$\mathbf{(viii)}$ If a system of equations over an $A-$algebra $B$ has a solution in each ring $B_{\mathfrak{m}}$ with $\mathfrak{m}$ a maximal ideal of $A$, then that system has a solution in the ring $B$. \\
$\mathbf{(ix)}$ If $\mathfrak{m}$ and $\mathfrak{m}'$ are distinct maximal ideals of $A$, then there exists an idempotent $e\in A$ such that $e\in\mathfrak{m}$ and $1-e\in\mathfrak{m}'$. \\
$\mathbf{(x)}$ If $\mathfrak{m}$ and $\mathfrak{m}'$ are distinct maximal ideals of $A$, then there exists an idempotent $e\in A$ such that $e\in\Ker\pi_{\mathfrak{m}}$ and $1-e\in\Ker\pi_{\mathfrak{m}'}$. \\
$\mathbf{(xi)}$ The connected components of $\Spec(A)$ are precisely of the form $\Lambda(\mathfrak{m})$ where $\mathfrak{m}$ is a maximal ideal of $A$. \\
$\mathbf{(xii)}$ The map $\lambda:\Max(A)\rightarrow\Sp(A)$ given by $\mathfrak{m}\rightsquigarrow(f\in\mathfrak{m}: f=f^{2})$ is a homeomorphism. \\
\end{theorem}

{\bf Proof.} In order to see the equivalences of $\mathbf{(i)}-\mathbf{(v)}$ consider \cite[Theorem 1.7]{McGovern}), and for the equivalences of $\mathbf{(i)}$, $\mathbf{(vi)}$ and $\mathbf{(vii)}$ see
\cite[Theorem 1.1]{Couchot}. \\
$\mathbf{(vi)}\Rightarrow\mathbf{(viii)}:$ Consider the system of equations
$f_{i}(x_{1},...,x_{n})=0$ with $i=1,...,d$ where $f_{i}(x_{1},...,x_{n})\in B[x_{1},...,x_{n}]$. (If $\phi:A\rightarrow B$ is the structure morphism then as usual $a.1_{B}=\phi(a)$ is simply denoted by $a$ for all $a\in A$). Using the calculus of fractions, then we may find a positive integer $N$ and polynomials $g_{i}(y_{1},...,y_{n},z_{1},...,z_{n})$ over $B$ such that $$f_{i}(b_{1}/s_{1},...,b_{n}/s_{n})=
g_{i}(b_{1},...,b_{n},s_{1},...,s_{n})/(s_{1}...s_{n})^{N}$$ for all $i$ and for every $b_{1},...,b_{n}\in B$ and $s_{1},...,s_{n}\in S$ where $S$ is a multiplicative subset of $A$. If the above system has a solution in each ring $B_{\mathfrak{m}}$, then there exist $b_{1},...,b_{n}\in B$ and $c,s_{1},...,s_{n}\in A\setminus\mathfrak{m}$ such that $$cg_{i}(b_{1},...,b_{n},s_{1},...,s_{n})=0$$ for all $i$. This leads us to consider $\mathscr{C}$ the collection of those opens $W$ of $\Max(A)$ such that there exists $b_{1},...,b_{n}\in B$ and $c,s_{1},...,s_{n}\in A\setminus(\bigcup\limits_{\mathfrak{m}\in W}\mathfrak{m})$ so that $$cg_{i}(b_{1},...,b_{n},s_{1},...,s_{n})=0.$$ Clearly $\mathscr{C}$ covers $\Max(A)$ and if $W\in\mathscr{C}$ then every open subset of $W$ is also a member of $\mathscr{C}$. Thus by Theorem \ref{Theorem IX}, we may find a finite number $W_{1},...,W_{q}\in\mathscr{C}$ of pairwise disjoint clopens of $\Max(A)$ such that $\Max(A)=\bigcup\limits_{k=1}^{q}W_{k}$. Using the retraction map $\gamma:\Spec(A)\rightarrow\Max(A)$ and Theorem \ref{Theorem II1}, then the map $f\rightsquigarrow D(f)\cap\Max(A)$ is a bijection from the set of idempotents of $A$ onto the set of clopens of $\Max(A)$. Therefore there exist orthogonal idempotents $e_{1},...,e_{q}\in A$ such that $W_{k}=D(e_{k})\cap\Max(A)$. Clearly $\sum\limits_{k=1}^{q}e_{k}$ is an idempotent and $D(\sum\limits_{k=1}^{q}e_{k})=\Spec(A)$. It follows that $\sum\limits_{k=1}^{q}e_{k}=1$. For each $k=1,...,q$ there exist $b_{1k},...,b_{nk}\in B$ and $c_{k}, s_{1k},...,s_{nk}\in A\setminus(\bigcup\limits_{\mathfrak{m}\in W_{k}}\mathfrak{m})$ such that
$$c_{k}g_{i}(b_{1k},...,b_{nk},s_{1k},...,s_{nk})=0$$ for all $i$. For each $j=1,...,n$ setting $b'_{j}:=\sum\limits_{k=1}^{q}e_{k}b_{jk}$ and $s'_{j}:=\sum\limits_{k=1}^{q}e_{k}s_{jk}$. Note that for each natural number $p\geqslant0$ we have then $(b'_{j})^{p}=\sum\limits_{k=1}^{q}e_{k}(b_{jk})^{p}$ and $(s'_{j})^{p}=\sum\limits_{k=1}^{q}e_{k}(s_{jk})^{p}$. We claim that $$c'g_{i}(b'_{1},...,b'_{n},s'_{1},...,s'_{n})=0$$ for all $i$ where $c':=\sum\limits_{k=1}^{q}e_{k}c_{k}$. Because fix $i$ and let $$g_{i}(y_{1},...,y_{n},z_{1},...,z_{n})=\sum\limits_{0\leqslant i_{1},....,i_{2n}<\infty}r_{i_{1},...,i_{2n}}y^{i_{1}}_{1}...
y^{i_{n}}_{n}z^{i_{n+1}}_{1}...z^{i_{2n}}_{n}.$$ Then $$c'g_{i}(b'_{1},...,b'_{n},s'_{1},...,s'_{n})=$$$$
c'\Big(\sum\limits_{0\leqslant i_{1},....,i_{2n}<\infty}r_{i_{1},...,i_{2n}}\big(\sum\limits_{k=1}^{q}
e_{k}(b_{1k})^{i_{1}}...
(b_{nk})^{i_{n}}(s_{1k})^{i_{n+1}}...(s_{nk})^{i_{2n}}\big)\Big)=$$$$
(\sum\limits_{t=1}^{q}e_{t}c_{t})\big(\sum\limits_{k=1}^{q}e_{k}
g_{i}(b_{1k},...,b_{nk},s_{1k},...,s_{nk})\big)=
$$$$\sum\limits_{k=1}^{q}e_{k}c_{k}
g_{i}(b_{1k},...,b_{nk},s_{1k},...,s_{nk})=0.$$
This establishes the claim. But $c'$ is invertible in $A$ since $c'\notin\mathfrak{m}$ for all $\mathfrak{m}\in\Max(A)$. Hence, $g_{i}(b'_{1},...,b'_{n},s'_{1},...,s'_{n})=0$ for all $i$. Similarly, each $s'_{j}$ is invertible in $A$. Therefore
$$f_{i}(b'_{1}/s'_{1},...,b'_{n}/s'_{n})=
 g_{i}(b'_{1},...,b'_{n},s'_{1},...,s'_{n})/(s'_{1}...s'_{n})^{N}=0$$ for all $i$.
Hence, the n-tuple $(b''_{1},...,b''_{n})\in B^{n}$ is a solution of the above system where $b_{j}'':=b'_{j}\phi(s'_{j})^{-1}$. \\
$\mathbf{(viii)}\Rightarrow\mathbf{(ii)}:$ It suffices to show that the system of equations: \[\left\{\
\begin{array}{ll}
X=X^{2}&\\
X=fY\\
1-X=(1-f)Z&\\
\end{array}\right. \]
has a solution in $A$. If $A$ is a local ring with the maximal ideal $\mathfrak{m}$, then the system having the solution $\big(0,0,1/(1-f)\big)$ or $(1,1/f,0)$, according as $f\in\mathfrak{m}$ or $f\notin\mathfrak{m}$. Using this, then by the hypothesis the system has a solution for every ring $A$ (not necessarily local). \\
$\mathbf{(vii)}\Rightarrow\mathbf{(ix)}:$ If $\mathfrak{m}$ and $\mathfrak{m}'$ are distinct maximal ideals of $A$, then by Theorem \ref{Theorem III}, $\Ker\pi_{\mathfrak{m}}+\Ker\pi_{\mathfrak{m}'}=A$. Thus there exists an idempotent $e\in\Ker\pi_{\mathfrak{m}}\subseteq\mathfrak{m}$ such that $e\notin\mathfrak{m}'$. It follows that $1-e\in\mathfrak{m}'$. \\
$\mathbf{(ix)}\Rightarrow\mathbf{(vi)}:$ It is clear. \\
$\mathbf{(ix)}\Leftrightarrow\mathbf{(x)}:$ Easy. \\
$\mathbf{(vi)}\Rightarrow\mathbf{(xi)}:$ If $\mathfrak{m}$ is a maximal ideal of $A$, then $A/\Ker\pi_{\mathfrak{m}}$ has no nontrivial idempotents since by Theorem \ref{Theorem III} (vi), it is canonically isomorphic to $A_{\mathfrak{m}}$. By (vii), $\Ker\pi_{\mathfrak{m}}$ is a regular ideal. It follows that $\Ker\pi_{\mathfrak{m}}$ is a max-regular ideal of $A$. Hence, by \cite[Theorem 3.17]{Abolfazl}, $V(\Ker\pi_{\mathfrak{m}})$ is a connected component of $\Spec(A)$. Conversely, if $C$ is a connected component of $\Spec(A)$, then $\gamma(C)$ is a connected subset of $\Max(A)$ where $\gamma:\Spec(A)\rightarrow\Max(A)$ is the retraction map which sends each prime ideal $\mathfrak{p}$ of $A$ into the unique maximal ideal of $A$ containing $\mathfrak{p}$. Therefore there exists a maximal ideal $\mathfrak{m}$ of $A$ such that $\gamma(C)=\{\mathfrak{m}\}$, because $\Max(A)$ is totally disconnected. We have $C\subseteq\gamma^{-1}(\{\mathfrak{m}\})=\Lambda(\mathfrak{m})=
V(\Ker\pi_{\mathfrak{m}})$. It follows that $C=V(\Ker\pi_{\mathfrak{m}})$. \\
$\mathbf{(xi)}\Rightarrow\mathbf{(vi)}:$ Clearly $A$ is a Gelfand ring, because distinct connected components are disjoint. The map $\phi:\Max(A)\rightarrow\pi_{0}\big(\Spec(A)\big)$ given by $\mathfrak{m}\rightsquigarrow\Lambda(\mathfrak{m})$ is bijective. It is also continuous and closed map, because $\phi=\pi\circ i$ and $\pi_{0}\big(\Spec(A)\big)$ is Hausdorff where $i:\Max(A)\rightarrow\Spec(A)$ is the canonical injection and $\pi:\Spec(A)\rightarrow\pi_{0}\big(\Spec(A)\big)$ is the canonical projection. Therefore, $\Max(A)$ is totally disconnected. \\
$\mathbf{(vii)}\Rightarrow\mathbf{(xii)}:$ The map $\lambda$ is continuous. It is also a closed map, because for any ring $A$, $\Max(A)$ is quasi-compact and $\Sp(A)$ is Hausdorff. If $\mathfrak{m}$ is a maximal ideal of $A$, then by the hypothesis, $\lambda(\mathfrak{m})=\Ker\pi_{\mathfrak{m}}$. This yields that $\lambda$ is injective, because $A$ is a Gelfand ring. If $M$ is a max-regular ideal of $A$, then there exists a maximal ideal $\mathfrak{m}$ of $A$ such that $M\subseteq\mathfrak{m}$. It follows that $M\subseteq\Ker\pi_{\mathfrak{m}}=\lambda(\mathfrak{m})$. This yields that $M=\lambda(\mathfrak{m})$. Hence, $\lambda$ is surjective. \\
$\mathbf{(xii)}\Rightarrow\mathbf{(vi)}:$ By the hypothesis, $\Max(A)$ is totally disconnected and it is the retraction of $\Spec(A)$. Hence, $A$ is a Gelfand ring. $\Box$ \\

\begin{corollary} The max-regular ideals of a clean ring $A$ are precisely of the form $\Ker\pi_{\mathfrak{m}}$ where $\mathfrak{m}$ is a maximal ideal of $A$. \\
\end{corollary}

{\bf Proof.} If $M$ is a max-regular ideal of $A$, then $V(M)$ is a connected component of $\Spec(A)$. Thus by Theorem \ref{Theorem VII} (xi), there exists a maximal ideal $\mathfrak{m}$ of $A$ such that $V(M)=\Lambda(\mathfrak{m})$. It follows that $M\subseteq\mathfrak{m}$. This yields that $M\subseteq\Ker\pi_{\mathfrak{m}}$. But $\Ker\pi_{\mathfrak{m}}$ is a regular ideal of $A$, see Theorem \ref{Theorem VII} (vii). Therefore $M=\Ker\pi_{\mathfrak{m}}$. By a similar argument, it is proven that if $\mathfrak{m}$ is a maximal ideal of $A$, then $\Ker\pi_{\mathfrak{m}}$ is a max-regular ideal of $A$. $\Box$ \\

\begin{corollary}\cite[Theorem 9]{Anderson}\label{Proposition V} If $A/\mathfrak{N}$ is a clean ring, then $A$ is a clean ring. \\
\end{corollary}

{\bf Proof.} It implies from Theorem \ref{Theorem VII} (xi). $\Box$ \\

\begin{corollary}\cite[Corollary 11]{Anderson}\label{Corollary VI} Every zero dimensional ring is a clean ring. \\
\end{corollary}

{\bf Proof.} If $A$ is a zero dimensional ring then by Theorem \ref{Theorem II} (iii), the Zariski and patch topologies over $\Max(A)$ are the same and so it is totally disconnected. Therefore by Theorem \ref{Theorem VII} (vi), $A$ is a clean ring. $\Box$ \\

\section{Mp-rings}

In this section mp-rings are characterized. \\

\begin{remark}\label{Remark VI} We observed that if $A$ is a Gelfand ring, then $\Max(A)$ is Zariski Hausdorff. Dually, if $A$ is a mp-ring then $\Min(A)$ is flat Hausdorff, because if $\mathfrak{p}$ and $\mathfrak{q}$ are distinct minimal primes of $A$ then $\mathfrak{p}+\mathfrak{q}=A$, thus there are $f\in\mathfrak{p}$ and $g\in\mathfrak{q}$ such that $f+g=1$ and so $V(f)\cap V(g)=\emptyset$. \\
\end{remark}

The following result is the culmination of mp-rings. \\

\begin{theorem}\label{Theorem V} For a ring $A$ the following statements are equivalent. \\
$\mathbf{(i)}$  $A$ is a mp-ring. \\
$\mathbf{(ii)}$ If $\mathfrak{p}$ and $\mathfrak{q}$ are distinct minimal prime ideals of $A$, then $\mathfrak{p}+\mathfrak{q}=A$. \\
$\mathbf{(iii)}$ $A/\mathfrak{N}$ is a mp-ring. \\
$\mathbf{(iv)}$ If $\mathfrak{p}$ is a minimal prime ideal of $A$, then $[\mathfrak{p}]=V(\mathfrak{p})$. \\
$\mathbf{(v)}$ $\Min(A)$ is the flat retraction of $\Spec(A)$. \\
$\mathbf{(vi)}$ $\Spec(A)$ is a normal space with respect to the flat topology.  \\
$\mathbf{(vii)}$ If $\mathfrak{p}$ is a minimal prime ideal of $A$, then $V(\mathfrak{p})$ is a flat closed subset of $\Spec(A)$. \\
$\mathbf{(viii)}$ The map $\eta:\Min(A)\rightarrow\Spec(A)/\sim_{R}$ given by $\mathfrak{p}\rightsquigarrow[\mathfrak{p}]$ is a homeomorphism. \\
$\mathbf{(ix)}$ If $fg=0$ for some $f,g\in A$, then there exists a natural number $n\geqslant1$ such that $\Ann(f^{n})+\Ann(g^{n})=A$. \\
$\mathbf{(x)}$ Every minimal prime ideal of $A$ is the radical of a unique pure ideal of $A$. \\
\end{theorem}

{\bf Proof.} The implications $\mathbf{(i)}\Rightarrow\mathbf{(ii)}\Rightarrow\mathbf{(iii)}
\Rightarrow\mathbf{(i)}$ are easy. \\
$\mathbf{(i)}\Rightarrow\mathbf{(iv)}:$ Let $\mathfrak{q}\in[\mathfrak{p}]$. There exists a minimal prime ideal $\mathfrak{p}'$ of $A$ such that $\mathfrak{p}'\subseteq\mathfrak{q}$. It follows that $\mathfrak{p}\sim_{R}\mathfrak{p}'$. Then by Remark \ref{Remark VIII}, there exists a finite set $\{\mathfrak{q}_{1},...,\mathfrak{q}_{n}\}$ of prime ideals of $A$ with $n\geqslant2$ such that $\mathfrak{q}_{1}=\mathfrak{p}$, $\mathfrak{q}_{n}=\mathfrak{p}'$ and $\mathfrak{q}_{i}+\mathfrak{q}_{i+1}\neq A$ for all $1\leqslant i\leqslant n-1$. By induction on $n$ we shall prove that $\mathfrak{p}=\mathfrak{p}'$. If $n=2$ then $\mathfrak{p}+\mathfrak{p}'\neq A$ and so by the hypothesis, $\mathfrak{p}=\mathfrak{p}'$. Assume that $n>2$. There exists a minimal prime ideal $\mathfrak{p}''$ of $A$ such that $\mathfrak{p}''\subseteq\mathfrak{q}_{n-1}$.
We have $\mathfrak{q}_{n-1}+\mathfrak{p}'\neq A$. Thus by the hypothesis, $\mathfrak{p}'=\mathfrak{p}''$. It follows that $\mathfrak{q}_{n-2}+\mathfrak{p}'\neq A$. Thus in the equivalence $\mathfrak{p}\sim_{R}\mathfrak{p}'$ the number of the involved primes is reduced to $n-1$. Therefore by the induction hypothesis, $\mathfrak{p}=\mathfrak{p}'$. \\
$\mathbf{(iv)}\Rightarrow\mathbf{(i)}:$ Let $\mathfrak{q}$ be a prime ideal of $A$ such that $\mathfrak{p}\subseteq\mathfrak{q}$ and $\mathfrak{p}'\subseteq\mathfrak{q}$ for some minimal primes $\mathfrak{p}$ and $\mathfrak{p}'$ of $A$. This implies that $\mathfrak{p}\sim_{R}\mathfrak{p}'$ and so $[\mathfrak{p}]=[\mathfrak{p}']$. This yields that $\mathfrak{p}=\mathfrak{p}'$. \\
$\mathbf{(i)}\Rightarrow\mathbf{(v)}:$ Consider the function $\gamma:\Spec(A)\rightarrow\Min(A)$ which sends each prime ideal $\mathfrak{p}$ of $A$ into the unique minimal prime ideal of $A$ which is contained in $\mathfrak{p}$. It suffices to show that $\gamma^{-1}\big(D(f)\cap\Min(A)\big)$ is a flat closed subset of $\Spec(A)$ for all $f\in A$. Let $E:=D(f)\cap\Min(A)$.
We show that $\gamma^{-1}(E)=\Ima\pi^{\ast}$ where $\pi:A\rightarrow S^{-1}A$ is the canonical ring map and $S=A\setminus\bigcup\limits_{f\notin\gamma(\mathfrak{p})}\mathfrak{p}$. Clearly $\gamma^{-1}(E)\subseteq\Ima\pi^{\ast}$. Conversely, if $\mathfrak{q}\in\Ima\pi^{\ast}$ then $\mathfrak{q}\subseteq\bigcup\limits_{f\notin\gamma(\mathfrak{p})}
\mathfrak{p}$. If $f\in\gamma(\mathfrak{q})$, then for each $\mathfrak{p}\in E$ there exist $x_{\mathfrak{p}}\in\mathfrak{p}$ and $y_{\mathfrak{p}}\in\gamma(\mathfrak{q})$ such that $x_{\mathfrak{p}}+y_{\mathfrak{p}}=1$. It follows that
$E\subseteq\bigcup\limits_{\mathfrak{p}\in E}
V(x_{\mathfrak{p}})$. But $E$ is a flat closed subset of $\Min(A)$ and for any ring $A$ the subspace $\Min(A)$ is flat quasi-compact.
This yields that $E$ is flat quasi-compact. Hence $E\subseteq\bigcup\limits_{i=1}^{n}
V(x_{i})=V(x)$ where $x=\prod\limits_{i=1}^{n}x_{i}$ and
$x_{i}:=x_{\mathfrak{p}_{i}}$ for all $i$. Thus we may find some $y\in\gamma(\mathfrak{q})$ such that $x+y=1$. Hence there exists a prime ideal $\mathfrak{p}$ such that $f\notin\gamma(\mathfrak{p})$ and $y\in\mathfrak{p}$. This yields that $1=x+y\in\mathfrak{p}$, a contradiction. Therefore $f\notin\gamma(\mathfrak{q})$. \\
$\mathbf{(v)}\Rightarrow\mathbf{(i)}:$ Let $\mathfrak{p}$ be a minimal prime ideal of $A$ such that $\mathfrak{p}\subseteq\mathfrak{q}$ for some prime ideal $\mathfrak{q}$ of $A$. By the hypothesis there exists a retraction map $\phi:\Spec(A)\rightarrow\Min(A)$. Clearly $\mathfrak{p}\in\Lambda(\mathfrak{q})=\overline{\{\mathfrak{q}\}}$ and so $\mathfrak{p}=\phi(\mathfrak{p})\in
\overline{\{\phi(\mathfrak{q})\}}=
\Lambda\big(\phi(\mathfrak{q})\big)\cap\Min(A)=
\{\phi(\mathfrak{q})\}$. It follows that $\mathfrak{p}=\phi(\mathfrak{q})$. \\
$\mathbf{(i)}\Rightarrow\mathbf{(vi)}:$ If $E$ and $F$ are disjoint flat closed subsets of $\Spec(A)$, then $\gamma(E)\cap\gamma(F)=\emptyset$, because flat closed subsets are stable under the generalization where $\gamma:\Spec(A)\rightarrow\Min(A)$ is the retraction map. \\
By (v), the map $\gamma$ is continuous. The space $\Min(A)$ is also flat Hausdorff, see Remark \ref{Remark VI}. Therefore $\gamma$ is a closed map. But $\Min(A)$ is a normal space, because it is well known that every compact space is a normal space. Thus there exist opens $\gamma(E)\subseteq U$ and $\gamma(F)\subseteq V$ in $\Min(A)$ such that $U\cap V=\emptyset$. It follows that $E\subseteq\gamma^{-1}(U)$ and $F\subseteq\gamma^{-1}(V)$ are disjoint opens of $\Spec(A)$. \\
$\mathbf{(vi)}\Rightarrow\mathbf{(ii)}:$ If $\mathfrak{p}$ and $\mathfrak{q}$ are distinct minimal prime ideals of $A$, then by the hypothesis there exists (finitely generated) ideals $I$ and $J$ of $A$ such that $\mathfrak{p}\in V(I)$, $\mathfrak{q}\in V(J)$ and $V(I)\cap V(J)=\emptyset$. This yields that $\mathfrak{p}+\mathfrak{q}=A$. \\
$\mathbf{(i)}\Rightarrow\mathbf{(vii)}:$ By \cite[Theorem 3.11]{Abolfazl}, it suffices to show that $V(\mathfrak{p})$ is stable under the generalization. Let $\mathfrak{q}'\subseteq\mathfrak{q}$ be prime ideals of $A$ such that $\mathfrak{p}\subseteq\mathfrak{q}$. There exists a minimal prime ideal $\mathfrak{p}'$ of $A$ such that $\mathfrak{p}'\subseteq\mathfrak{q}'$. This yields that $\mathfrak{p}=\mathfrak{p}'$, since $A$ is a mp-ring. \\
$\mathbf{(vii)}\Rightarrow\mathbf{(ii)}:$ If $\mathfrak{p}+\mathfrak{q}\neq A$, then there exists a maximal ideal $\mathfrak{m}$ of $A$ such that $\mathfrak{p}+\mathfrak{q}\subseteq\mathfrak{m}$. It follows that $\mathfrak{m}\in V(\mathfrak{p})$. This yields that $\mathfrak{q}\in V(\mathfrak{p})$, because by the hypothesis, $V(\mathfrak{p})$ is stable under the generalization. This implies that $\mathfrak{p}=\mathfrak{q}$. But this is a contradiction. \\
$\mathbf{(i)}\Rightarrow\mathbf{(viii)}:$ This is proved exactly like the proof of the implication (i)$\Rightarrow$(xi) of Theorem \ref{Theorem III}. \\
$\mathbf{(viii)}\Rightarrow\mathbf{(i)}:$ It is proved exactly like the implication (iv)$\Rightarrow$(i). \\
$\mathbf{(i)}\Rightarrow\mathbf{(ix)}:$ Suppose $\Ann(f^{n})+\Ann(g^{n})\neq A$ for all $n\geqslant1$. Then $I:=\sum\limits_{n\geqslant1}\big(\Ann(f^{n})+\Ann(g^{n})\big)$ is a proper ideal of $A$. So there exists a maximal ideal $\mathfrak{m}$ of $A$
such that $I\subseteq\mathfrak{m}$. Let $\mathfrak{p}$ be the unique minimal prime ideal of $A$ contained in $\mathfrak{m}$. If $\mathfrak{N}$ is the nil-radical of $A_{\mathfrak{m}}$, then
Clearly $f/1\notin\mathfrak{N}=\mathfrak{p}A_{\mathfrak{m}}$, because $\Ann(f^{n})\subseteq\mathfrak{m}$ for all $n\geqslant1$. Therefore $f\notin\mathfrak{p}$. By a similar argument, we get that $g\notin\mathfrak{p}$. But this is a contradiction, since $fg=0$. \\
$\mathbf{(ix)}\Rightarrow\mathbf{(ii)}:$ By Lemma \ref{Lemma III},
there exist $f\in A\setminus\mathfrak{p}$ and $g\in A\setminus\mathfrak{q}$ such that $fg=0$. Thus by the hypothesis, $\Ann(f^{n})+\Ann(g^{n})=A$ for some $n\geqslant1$. Hence, there exist $x\in\Ann(f^{n})$ and $y\in\Ann(g^{n})$ such that $1=x+y$. It follows that $x\in\mathfrak{p}$ and $y\in\mathfrak{q}$ and so
$\mathfrak{p}+\mathfrak{q}=A$. \\
$\mathbf{(vii)}\Rightarrow\mathbf{(x)}:$ By \cite[Theorem 3.11]{Abolfazl}, $V(\mathfrak{p})$ is stable under the generalization. Thus by \cite[Theorem 3.2]{Abolfazl 2020}, there exists a unique pure ideal $I$ of $A$ such that $V(\mathfrak{p})=V(I)$. So $\mathfrak{p}=\sqrt{I}$. \\
$\mathbf{(x)}\Rightarrow\mathbf{(ii)}:$ There exists some $f\in\mathfrak{p}$ such that $f\notin\mathfrak{q}$. By the hypothesis, there exists a pure ideal $I$ of $A$ such that $\mathfrak{p}=\sqrt{I}$. So by Theorem \ref{Proposition I}, there exist some $g\in I$ and a natural number $n\geqslant1$ such that $f^{n}(1-g)=0$. It follows that $1-g\in\mathfrak{q}$. Thus $\mathfrak{p}+\mathfrak{q}=A$. $\Box$ \\

If $\mathfrak{p}$ is a minimal prime ideal of a mp-ring $A$, then the unique pure ideal in Theorem \ref{Theorem V} (x) associated to $\mathfrak{p}$, is the largest pure ideal of $A$ which is contained in $\Ker\pi_{\mathfrak{p}}$. \\
We give a new proof to the following technical result. Our proof essentially differs with \cite[\S4]{Couchot} from the scratch. \\

\begin{theorem}\label{Theorem Tariz-Couch} Let $A$ be a mp-ring. Then the retraction map $\gamma:\Spec(A)\rightarrow\Min(A)$ is Zariski continuous if and only if $\Min(A)$ is Zariski compact. \\
\end{theorem}

{\bf Proof.} If $\gamma$ is Zariski continuous, then clearly $\Min(A)$ is Zariski quasi-compact. For any ring $A$, then $\Min(A)$ is Zariski Hausdorff, see Lemma \ref{Lemma III}. Thus $\Min(A)$ is Zariski compact. To see the converse it suffices to show that $\gamma^{-1}(U)$ is a Zariski open of $\Spec(A)$ where $U=\Min(A)\cap D(f)$ and $f\in A$. It is well known and easy to see that for any ring $A$, then $U$ is a Zarisiki clopen of $\Min(A)$, see \cite[Proposition 4.2]{Tarizadeh 303}. Thus there exists an ideal $I$ of $A$ such that $U=\Min(A)\cap V(I)=\bigcap\limits_{g\in I}\Min(A)\cap V(g)$. It follows that $U^{c}=\Min(A)\setminus U=\bigcup\limits_{g\in I}\Min(A)\cap D(g)$. By the hypothesis, $\Min(A)$ is Zariski quasi-compact. Thus there exist a finitely many elements $g_{1},...,g_{n}\in I$ such that $U^{c}=\bigcup\limits_{i=1}^{n}\Min(A)\cap D(g_{i})$, since every closed subset of a quasi-compact space is quasi-compact. This yields that $U=\Min(A)\cap V(J)$ where $J=(g_{1},...,g_{n})$. Thus $U$ is a flat open of $\Min(A)$. It is also a flat closed subset of $\Min(A)$. Therefore $\gamma^{-1}(U)$ is a flat clopen of $\Spec(A)$, because by Theorem \ref{Theorem V} (v), the map $\gamma$ is flat continuous. But Zariski clopens and flat clopens of $\Spec(A)$ are the same, see \cite[Corollary 3.12]{Abolfazl}. Hence, $\gamma^{-1}(U)$ is a Zariski open. $\Box$ \\

If $A$ is a ring, then clearly $\Ann(f)+\Ann(g)\subseteq\Ann(fg)$ for all $f,g\in A$. In Theorem \ref{Theorem VI}, it is shown that the equality holds if and only if $A$ is a reduced mp-ring. In Theorem \ref{Theorem VI}, we also improve \cite[Theorem 2.6]{Al Ezeh 2},
\cite[Lemma $\beta$ ]{Artico-Marconi} and \cite[Theorem 2.3]{Bhattacharjee} by adding the condition (viii) as a new equivalent. \\

\begin{theorem}\label{Theorem VI} For a ring $A$ the following statements are equivalent. \\
$\mathbf{(i)}$ $A$ is a reduced mp-ring. \\
$\mathbf{(ii)}$ If $fg=0$ for some $f,g\in A$, then $\Ann(f)+\Ann(g)=A$. \\
$\mathbf{(iii)}$  $\Ann(f)+\Ann(g)=\Ann(fg)$ for all $f,g\in A$. \\
$\mathbf{(iv)}$ $\Ann(f)$ is a pure ideal for all $f\in A$. \\
$\mathbf{(v)}$ Every principal ideal of $A$ is a flat $A-$module. \\
$\mathbf{(vi)}$ If $\mathfrak{p}$ is a prime ideal of $A$, then $A_{\mathfrak{p}}$ is an integral domain. \\
$\mathbf{(vii)}$ If $\mathfrak{m}$ is a maximal ideal of $A$, then $A_{\mathfrak{m}}$ is an integral domain. \\
$\mathbf{(viii)}$ Every minimal prime ideal of $A$ is a pure ideal. \\
\end{theorem}

{\bf Proof.} $\mathbf{(i)}\Rightarrow\mathbf{(ii)}:$ By Theorem \ref{Theorem V} (ix), $\Ann(f^{n})+\Ann(g^{n})=A$ for some $n\geqslant1$. But $\Ann(f)=\Ann(f^{k})$ for all $f\in A$ and all $k\geqslant1$, since $A$ is reduced. \\
$\mathbf{(ii)}\Rightarrow\mathbf{(i)}:$ By Theorem \ref{Theorem V} (ix), $A$ is a mp-ring. Let $f$ be a nilpotent element of $A$. Thus there exists the least positive natural number $n$ such that $f^{n}=0$. We show that $n=1$. If $n>1$ then by the hypothesis, $\Ann(f^{n-1})=\Ann(f)+\Ann(f^{n-1})=A$. It follows that $f^{n-1}=0$. But this is in contradiction with the minimality of $n$. \\
$\mathbf{(ii)}\Rightarrow\mathbf{(iii)}:$ If $a\in\Ann(fg)$ then $(af)g=0$. Thus by the hypothesis, $\Ann(af)+\Ann(g)=A$. Hence there are $b\in\Ann(af)$ and $c\in\Ann(g)$ such that $b+c=1$. We have $a=ab+ac$, $ab\in\Ann(f)$ and $ac\in\Ann(g)$. Thus $a\in\Ann(f)+\Ann(g)$. The implications $\mathbf{(iii)}\Rightarrow\mathbf{(ii)}$,
$\mathbf{(ii)}\Leftrightarrow\mathbf{(iv)}\Leftrightarrow\mathbf{(v)}$ and $\mathbf{(vi)}\Rightarrow\mathbf{(vii)}\Rightarrow\mathbf{(ii)}$ are easy. \\
$\mathbf{(i)}\Rightarrow\mathbf{(vi)}:$ If $\mathfrak{q}$ is the minimal prime ideal of $A$ contained in $\mathfrak{p}$, then $\mathfrak{q}A_{\mathfrak{p}}=0$. Hence, $A_{\mathfrak{p}}$ is a domain. \\
$\mathbf{(i)}\Rightarrow\mathbf{(viii)}:$ Let $\mathfrak{p}$ be a minimal prime ideal of $A$ and $f\in\mathfrak{p}$. If $\Ann(f)+\mathfrak{p}\neq A$, then there exists a maximal ideal $\mathfrak{m}$ of $A$ such that $\Ann(f)+\mathfrak{p}\subseteq\mathfrak{m}$. By the hypotheses, $\mathfrak{p}A_{\mathfrak{m}}=0$. Hence there exists some $g\in A\setminus\mathfrak{m}$ such that $fg=0$. But this is a contradiction. Thus by Theorem \ref{Proposition I}, $\mathfrak{p}$ is a pure ideal. \\
$\mathbf{(viii)}\Rightarrow\mathbf{(i)}:$ By Theorem \ref{Theorem V} (x), $A$ is a mp-ring. If $f\in A$ is nilpotent, then $\Ann(f)=A$. If not, then there exists a maximal ideal $\mathfrak{m}$ of $A$ such that $\Ann(f)\subseteq\mathfrak{m}$. There exists a minimal prime ideal $\mathfrak{p}$ of $A$ such that $\mathfrak{p}\subseteq\mathfrak{m}$. Clearly $f\in\mathfrak{p}$ and so by the hypothesis, there is some $g\in\mathfrak{p}$ such that $1-g\in\Ann(f)$ which is a contradiction. Thus $f=0$. Hence, $A$ is reduced. $\Box$ \\

In the literature, a ring $A$ is said to be a p.f. ring if every principal ideal of $A$ is a flat $A-$module. By Theorem \ref{Theorem VI}, p.f. rings and reduced mp-rings are the same. \\

\begin{corollary}\label{Corollary IX} If $\mathfrak{p}$ is a prime ideal of a reduced mp-ring $A$, then $\Ker\pi_{\mathfrak{p}}$ is the unique minimal prime ideal of $A$ which is contained in $\mathfrak{p}$.  \\
\end{corollary}

{\bf Proof.} Clearly $\Ker\pi_{\mathfrak{p}}\subseteq\mathfrak{p}$, and by Theorem \ref{Theorem VI} (ii), it is a prime ideal. Its minimality is easily verified and uniqueness is obvious, since $A$ is a mp-ring. $\Box$ \\

\section{Pure ideals of reduced Gelfand rings and mp-rings}

In this section, new characterizations for pure ideals of reduced Gelfand rings and mp-rings are given. \\
The following result generalizes \cite[Theorem 1.8]{Al Ezeh} to any ring. \\

\begin{lemma}\label{Lemma I23} If $I$ is a pure ideal of a ring $A$, then: $$I=\bigcap\limits_{\mathfrak{m}\in\Max(A)\cap V(I)}\Ker\pi_{\mathfrak{m}}.$$ \\
\end{lemma}

{\bf Proof.} If $\mathfrak{m}$ is a maximal ideal of $A$ containing $I$, then $I\subseteq\Ker\pi_{\mathfrak{m}}$, because if $f\in I$ then by Theorem \ref{Proposition I}, there exists some $g\in I$ such that $(1-g)f=0$, clearly $1-g\in A\setminus\mathfrak{m}$ and so $f\in\Ker\pi_{\mathfrak{m}}$. Conversely, if $f\in\bigcap\limits_{\mathfrak{m}\in\Max(A)\cap V(I)}\Ker\pi_{\mathfrak{m}}$ then to prove $f\in I$ it suffices to show that $\Ann(f)+I=A$. If $\Ann(f)+I\neq A$, then there exists a maximal ideal $\mathfrak{m}'$ of $A$ such that $\Ann(f)+I\subseteq\mathfrak{m}'$ and so $f\notin\Ker\pi_{\mathfrak{m}'}$ which is a contradiction. $\Box$ \\

\begin{theorem}\label{Theorem 19} The pure ideals of a reduced Gelfand ring $A$ are precisely of the form $\bigcap\limits_{\mathfrak{m}\in\Max(A)\cap E}\Ker\pi_{\mathfrak{m}}$ where $E$ is a Zariski closed subset of $\Spec(A)$. \\
\end{theorem}

{\bf Proof.} It follows from Lemma \ref{Lemma I23} and \cite[Theorem 1.9]{Al Ezeh}. $\Box$ \\

\begin{theorem}\label{Theorem 20} The pure ideals of a reduced mp-ring $A$ are precisely of the form $\bigcap\limits_{\mathfrak{p}\in\Min(A)\cap E}\mathfrak{p}$ where $E$ is a flat closed subset of $\Spec(A)$. \\
\end{theorem}

{\bf Proof.} If $I$ is a pure ideal of $A$, then $V(I)$ is a flat closed subset of $\Spec(A)$ and we have $I=\bigcap\limits_{\mathfrak{p}\in\Min(A)\cap V(I)}\mathfrak{p}$. Because if we choose $f$ in the intersection, then to prove $f\in I$ it suffices to show that $\Ann(f)+I=A$. If $\Ann(f)+I\neq A$, then $\Ann(f)+I\subseteq\mathfrak{m}$ for some maximal ideal $\mathfrak{m}$ of $A$. There exists a minimal prime ideal $\mathfrak{q}$ of $A$ such that $\mathfrak{q}\subseteq\mathfrak{m}$. But $I\subseteq\mathfrak{q}$, since $I$ is a pure ideal and so by Theorem \ref{Proposition I}, $V(I)$ is stable under the generalization. Hence,  $f\in\mathfrak{q}$. By Theorem \ref{Theorem VI}, $\mathfrak{q}$ is a pure ideal. Thus by Theorem \ref{Proposition I}, there exists some $g\in\mathfrak{q}$ such that $(1-g)f=0$. It follows that $1=(1-g)+g\in\mathfrak{m}$ which is a contradiction. Conversely, setting $J:=\bigcap\limits_{\mathfrak{p}\in\Min(A)\cap E}\mathfrak{p}$ where $E$ is a flat closed subset of $\Spec(A)$. Let $f\in J$. Then for each $\mathfrak{p}\in\Min(A)\cap E$ there exists some $c_{\mathfrak{p}}\in\mathfrak{p}$ such that $f=fc_{\mathfrak{p}}$, because by Theorem \ref{Theorem VI}, $\mathfrak{p}$ is a pure ideal of $A$. We claim that $J+(1-c_{\mathfrak{p}}: \mathfrak{p}\in\Min(A)\cap E)$ is the unit ideal of $A$. If not, then it is contained in a maximal ideal $\mathfrak{m}$ of $A$. There exists a minimal prime ideal $\mathfrak{q}$ of $A$ such that $\mathfrak{q}\subseteq\mathfrak{m}$. Clearly $\mathfrak{q}\notin E$. Thus for each $\mathfrak{p}\in\Min(A)\cap E$ there exist $x_{\mathfrak{p}}\in\mathfrak{p}$ and $y_{\mathfrak{p}}\in\mathfrak{q}$ such that $x_{\mathfrak{p}}+y_{\mathfrak{p}}=1$, because $A$ is a mp-ring. We have $E\subseteq\bigcup\limits_{\mathfrak{p}\in\Min(A)\cap E}V(x_{\mathfrak{p}})$, because every flat closed subset of $\Spec(A)$ is stable under the generalization. But every closed subset of a quasi-compact space is quasi-compact. Therefore there exist a finite number $\mathfrak{p}_{1},...,\mathfrak{p}_{n}\in\Min(A)\cap E$ such that $E\subseteq\bigcup\limits_{i=1}^{n}V(x_{i})$ where $x_{i}:=x_{\mathfrak{p}_{i}}$ for all $i$. There exists some $y\in\mathfrak{q}$ such that $x+y=1$ where $x=\prod\limits_{i=1}^{n}x_{i}$. Clearly $x\in J$. It follows that $1=x+y\in J+\mathfrak{q}\subseteq\mathfrak{m}$. This is a contradiction. Therefore $J+(1-c_{\mathfrak{p}}: \mathfrak{p}\in\Min(A)\cap E)=A$. Thus we may write $1=g+\sum\limits_{k}a_{k}(1-c_{k})$ where $g\in J$ and $c_{k}:=c_{\mathfrak{p}_{k}}$ for all $k$. It follows that $f=fg$. Therefore by Theorem \ref{Proposition I}, $J$ is a pure ideal. $\Box$ \\

\begin{corollary}\cite[Theorems 2.4 and 2.5]{Al Ezeh 2}\label{Corollary 29} The pure ideals of a reduced mp-ring $A$ are precisely of the form $\bigcap\limits_{\mathfrak{m}\in\Max(A)\cap V(I)}\Ker\pi_{\mathfrak{m}}$ where $I$ is an ideal of $A$. \\
\end{corollary}

{\bf Proof.} The implication ``$\Rightarrow$" implies from Lemma \ref{Lemma I23}. To prove the converse, first we claim that $\{\Ker\pi_{\mathfrak{m}}: \mathfrak{m}\in\Max(A)\cap V(I)\}=\{\mathfrak{p}: \mathfrak{p}\in\Min(A)\cap E\}$ where $E=\Ima\pi^{\ast}$ and $\pi:A\rightarrow S^{-1}A$ is the canonical ring map with $S=1+I$. The inclusion $\subseteq$ implies from Corollary \ref{Corollary IX}. To see the reverse inclusion, if $\mathfrak{p}\in\Min(A)\cap E$ then $I+\mathfrak{p}\neq A$. Thus there exists a maximal ideal $\mathfrak{m}$ of $A$ such that $I+\mathfrak{p}\subseteq\mathfrak{m}$. Then by Corollary \ref{Corollary IX}, $\mathfrak{p}=\Ker\pi_{\mathfrak{m}}$. Hence the claim is established. This yields that $\bigcap\limits_{\mathfrak{m}\in\Max(A)
\cap V(I)}\Ker\pi_{\mathfrak{m}}=\bigcap\limits_{\mathfrak{p}\in\Min(A)\cap E}\mathfrak{p}$. Now by Theorem \ref{Theorem 20}, the assertion is concluded. $\Box$ \\

\begin{corollary}\label{Corollary 1909} If $A$ is a reduced mp-ring, then $\Ann(f)\cap J_{f}=0$ for all $f\in A$ where $J_{f}:=\bigcap\limits_{\mathfrak{m}\in\Max(A)\cap V(f)}\Ker\pi_{\mathfrak{m}}$. \\
\end{corollary}

{\bf Proof.} By Theorem \ref{Theorem VI} and Corollary \ref{Corollary 29},  $\Ann(f)\cap J_{f}$ is a pure ideal and it is contained in the Jacobson radical of $A$. Hence, $\Ann(f)\cap J_{f}=0$. $\Box$ \\

The above corollary together with \cite[Lemma 3.4]{Al Ezeh 2} provide a short and straight proof for \cite[Theorem 3.5]{Al Ezeh 2}. \\

\begin{corollary}\label{Corollary 1910} If $\mathfrak{p}$ is a minimal prime ideal of a reduced mp-ring $A$, then $\mathfrak{p}=\sum\limits_{f\in\mathfrak{p}}J_{f}$. \\
\end{corollary}

{\bf Proof.} If $f\in\mathfrak{p}$ then there exists some $h\in\mathfrak{p}$ such that $f(1-h)=0$. This yields that $f\in J_{h}$. The reverse inclusion is deduced from Corollary \ref{Corollary 1909}. $\Box$ \\

\section{Purified rings}

The following definition introduces a new class of commutative rings. \\

\begin{definition} We call a ring $A$ a \emph{purified ring} if for all distinct minimal prime ideals $\mathfrak{p}$ and $\mathfrak{q}$ there exists an idempotent in $A$ which lies in $\mathfrak{p}$, but not in $\mathfrak{q}$. \\
\end{definition}

Obviously every purified ring is a mp-ring. The above definition, in the light of Theorem \ref{Theorem VII} (ix), is the dual notion of clean ring. \\
It is important to notice that, in order to get the dual notion of clean ring, the initial definition of clean ring can not be dualized by replacing ``product'' instead of ``sum''. Because each element of a ring can be written as a ``product'' of an invertible and an idempotent elements of that ring if and only if it is absolutely flat (von-Neumann regular) ring, see \cite[page 1]{Anderson}. This in particular yields that every absolutely flat ring is again a clean ring, see Corollary \ref{Corollary VI} or \cite[Theorem 10]{Anderson}. \\
Purified rings are stable under taking localizations. A finite product of rings is a purified ring if and only if each factor is a purified ring. \\

\begin{proposition} Every zero dimensional ring is a purified ring. \\
\end{proposition}

{\bf Proof.} It follows from Theorem \ref{Theorem VII} (ix) and Corollary \ref{Corollary VI}. $\Box$ \\

\begin{theorem}\label{Proposition IV089} Let $A$ be a ring. Then $A$ is a purified ring if and only if $A/\mathfrak{N}$ is a purified ring. \\
\end{theorem}

{\bf Proof.} Let $A/\mathfrak{N}$ be a purified ring and let $\mathfrak{p}$ and $\mathfrak{q}$ be distinct minimal prime ideals of $A$. Then there exists an idempotent $f+\mathfrak{N}\in A/\mathfrak{N}$ such that $f\in\mathfrak{p}$ and $1-f\in\mathfrak{q}$. Using Theorem \ref{Theorem II1}, then it is not hard to see that the idempotents of a ring $A$ can be lifted modulo its nil-radical.
So there exists an idempotent $e\in A$ such that $f-e\in\mathfrak{N}$. It follows that $e\in\mathfrak{p}$ and $1-e=(1-f)+(f-e)\in\mathfrak{q}$. $\Box$ \\

\begin{lemma}\label{Lemma nice 2} Let $I$ be a regular ideal of a ring $A$. If $f\in I$, then there exists an idempotent $e\in I$ such that $f=fe$. \\
\end{lemma}

{\bf Proof.} There exist a finitely many idempotents $e_{1},...,e_{n}\in I$ such that $f\in(e_{1},...,e_{n})$. But $e_{1}+e_{2}-e_{1}e_{2}$ is an idempotent and $(e_{1},e_{2})=(e_{1}+e_{2}-e_{1}e_{2})$. Thus by the induction, there exists an idempotent $e\in I$ such that $(e_{1},...,e_{n})=Ae$. Therefore $f=re$ for some $r\in A$. But we have $f(1-e)=re(1-e)=0$ and so $f=fe$. $\Box$ \\

The following result is the culmination of purified rings. \\

\begin{theorem}\label{Theorem X} For a reduced ring $A$ the following statements are equivalent. \\
$\mathbf{(i)}$ $A$ is a purified ring. \\
$\mathbf{(ii)}$ $A$ is a mp-ring and $\Min(A)$ is totally disconnected with respect to the flat topology. \\
$\mathbf{(iii)}$ Every minimal prime ideal of $A$ is a regular ideal. \\
$\mathbf{(iv)}$ The connected components of $\Spec(A)$ are precisely of the form $V(\mathfrak{p})$ where $\mathfrak{p}$ is a minimal prime ideal of $A$. \\
$\mathbf{(v)}$ If a system of equations over $A$ has a solution in each ring $A/\mathfrak{p}$ with $\mathfrak{p}$ a minimal prime of $A$, then that system has a solution in $A$. \\
$\mathbf{(vi)}$ The idempotents of $A$ can be lifted along each localization of $A$. \\
$\mathbf{(vii)}$ The collection of $V(e)\cap\Min(A)$ with $e\in A$ an idempotent forms a basis for the induced flat topology on $\Min(A)$. \\
$\mathbf{(viii)}$ $A$ is a mp-ring and every pure ideal of $A$ is a regular ideal. \\
$\mathbf{(ix)}$ The max-regular ideals of $A$ are precisely the minimal prime ideals of $A$. \\
$\mathbf{(x)}$ For each $f\in A$, $\Ann(f)$ is a regular ideal of $A$. \\
$\mathbf{(xi)}$ If $fg=0$ for some $f,g\in A$, then there exists an idempotent $e\in A$ such that $f=fe$ and $g=g(1-e)$. \\
\end{theorem}

{\bf Proof.} $\mathbf{(i)}\Rightarrow\mathbf{(ii)}:$ If $\mathfrak{p}$ and $\mathfrak{q}$ are distinct minimal primes of $A$, then there exists an idempotent $e\in A$ such that $\mathfrak{p}\in V(e)$ and $\mathfrak{q}\in V(1-e)$. We also have $V(e)\cup V(1-e)=\Spec(A)$. Therefore $\Min(A)$ is totally disconnected with respect to the flat topology. \\
$\mathbf{(ii)}\Rightarrow\mathbf{(iii)}:$ Let $\mathfrak{p}$ be a minimal prime of $A$ and $f\in\mathfrak{p}$. By Remark \ref{Remark VI}, $\Min(A)$ is flat Hausdorff. It is also flat quasi-compact. Therefore by Theorem \ref{Theorem IX}, there exists a clopen $U\subseteq\Min(A)$ such that $\mathfrak{p}\in U\subseteq V(f)\cap\Min(A)$. Then by Theorem \ref{Theorem II1}, there exists an idempotent $e\in A$ such that $\mathfrak{p}\in V(e)=\gamma^{-1}(U)$ where $\gamma:\Spec(A)\rightarrow\Min(A)$ is the retraction map, see Theorem \ref{Theorem V}. We have $\gamma^{-1}(U)\subseteq V(f)$. Thus there exist a natural number $n\geqslant1$ and an element $a\in A$ such that $f^{n}=ae$. It follows that $1-e\in\Ann(f^{n})$. But $\Ann(f^{n})=\Ann(f)$, since $A$ is reduced. Therefore $f=fe$. \\
$\mathbf{(iii)}\Rightarrow\mathbf{(i)}:$ Let $\mathfrak{p}$ and $\mathfrak{q}$ be distinct minimal primes of $A$. Then there exists an idempotent $e\in\mathfrak{p}$ such that $e\notin\mathfrak{q}$. It follows that $1-e\in\mathfrak{q}$. \\
$\mathbf{(ii)}\Rightarrow\mathbf{(iv)}:$ If $\mathfrak{p}$ is a minimal prime of $A$, then it is a max-regular ideal of $A$, see the implication (ii)$\Rightarrow$(iii). Thus by \cite[Theorem 3.17]{Abolfazl}, $V(\mathfrak{p})$ is a connected component of $A$. Conversely, if $C$ is a connected component of $\Spec(A)$, then there exists a minimal prime $\mathfrak{p}$ of $A$ such that $\gamma(C)=\{\mathfrak{p}\}$. But we have $C\subseteq\gamma^{-1}(\{\mathfrak{p}\})=V(\mathfrak{p})$. It follows that $C=V(\mathfrak{p})$. \\
$\mathbf{(iv)}\Rightarrow\mathbf{(ii)}:$ Clearly $A$ is a mp-ring, because distinct connected components are disjoint. The map $\Min(A)\rightarrow\pi_{0}\big(\Spec(A)\big)$ given by $\mathfrak{p}\rightsquigarrow V(\mathfrak{p})$ is a homeomorphism. Hence, $\Min(A)$ is totally disconnected with respect to the flat topology. \\
$\mathbf{(ii)}\Rightarrow\mathbf{(v)}:$ Assume that the system of equations $f_{i}(x_{1},...,x_{n})=0$ over $A$ has a solution in each ring $A/\mathfrak{p}$. Thus for each minimal prime ideal $\mathfrak{p}$ of $A$ there exist $b_{1},...,b_{n}\in A$ such that $f_{i}(b_{1},...,b_{n})\in\mathfrak{p}$ for all $i$. This leads us to consider $\mathscr{C}$, the collection of those opens $W$ of $\Min(A)$ such that there exist $b_{1},...,b_{n}\in A$ so that $f_{i}(b_{1},...,b_{n})\in\bigcap\limits_{\mathfrak{p}\in W}\mathfrak{p}$ for all $i$. Clearly $\mathscr{C}$ covers $\Min(A)$, and if $W\in\mathscr{C}$ then every open subset of $W$ is also a member of $\mathscr{C}$. Thus by Theorem \ref{Theorem IX}, we may find a finite number $W_{1},...,W_{q}\in\mathscr{C}$
of pairwise disjoint clopens of $\Min(A)$ such that $\Min(A)=\bigcup\limits_{k=1}^{q}W_{k}$. Using Theorem \ref{Theorem II1} and the retraction map $\gamma:\Spec(A)\rightarrow\Min(A)$ of Theorem \ref{Theorem V}, then the map $f\rightsquigarrow V(f)\cap\Min(A)$ is a bijection from the set of idempotents of $A$ onto the set of clopens of $\Min(A)$. Therefore there are orthogonal idempotents $e_{1},...,e_{q}\in A$ such that $W_{k}=V(1-e_{k})\cap\Min(A)$. Thus $\sum\limits_{k=1}^{q}e_{k}$ is an idempotent and $D(\sum\limits_{k=1}^{q}e_{k})=\Spec(A)$. It follows that $\sum\limits_{k=1}^{q}e_{k}=1$. For each $k=1,...,q$ there exist $b_{1k},...,b_{nk}\in A$ such that $f_{i}(b_{1k},...,b_{nk})\in\bigcap\limits_{\mathfrak{p}\in W_{k}}\mathfrak{p}$ for all $i$. For each $j=1,...,n$ setting $b'_{j}:=\sum\limits_{k=1}^{q}e_{k}b_{jk}$. Note that if $p\geqslant0$ is a natural number, then $(b'_{j})^{p}=\sum\limits_{k=1}^{q}e_{k}(b_{jk})^{p}$. It follows that $f_{i}(b'_{1},...,b'_{n})=\sum\limits_{k=1}^{q}e_{k}f_{i}(b_{1k},...,b_{nk})$ for all $i$. Now if $\mathfrak{p}$ is a minimal prime of $A$ then $\mathfrak{p}\in W_{t}$ for some $t$. We have $e_{t}f_{i}(b'_{1},...,b'_{n})=e_{t}f_{i}(b_{1t},...,b_{nt})\in\mathfrak{p}$.
This implies that $f_{i}(b'_{1},...,b'_{n})\in\mathfrak{p}$. Therefore $f_{i}(b'_{1},...,b'_{n})\in
\bigcap\limits_{\mathfrak{p}\in\Min(A)}\mathfrak{p}=0$ for all $i$. \\
$\mathbf{(v)}\Rightarrow\mathbf{(vi)}:$ Let $S$ be a multiplicative subset of $A$. If $a/s\in S^{-1}A$ is an idempotent, then there exists some $t\in S$ such that $ast(a-s)=0$.
It suffices to show that the following system: \[\left\{\
\begin{array}{ll}
X=X^{2}&\\
st(a-sX)=0\\
\end{array}\right. \]
has a solution in $A$. If $A$ is an integral domain, then the above system having the solution $0_{A}$ or $1_{A}$, according as $ast=0$ or $a=s$. Using this, then by the hypothesis the above system has a solution for every ring $A$ (not necessarily domain). \\
$\mathbf{(vi)}\Rightarrow\mathbf{(i)}:$ If $\mathfrak{p}$ and $\mathfrak{q}$ are distinct minimal prime ideals of $A$, then by the prime avoidance lemma, $\Spec(S^{-1}A)=\{S^{-1}\mathfrak{p}, S^{-1}\mathfrak{q}\}$ where $S=A\setminus(\mathfrak{p}\cup\mathfrak{q})$. Clearly $S^{-1}\mathfrak{p}$ and $S^{-1}\mathfrak{q}$ are distinct maximal ideals of $S^{-1}A$. Therefore by Theorem \ref{Theorem II1}, there exists an idempotent $f\in S^{-1}A$ such that $D(f)=\{S^{-1}\mathfrak{q}\}$ and $D(1-f)=\{S^{-1}\mathfrak{p}\}$. By the hypothesis, there exists an idempotent $e\in A$ such that $e/1=f$. It follows that $e\in\mathfrak{p}$ and $1-e\in\mathfrak{q}$. \\
$\mathbf{(ii)}\Rightarrow\mathbf{(vii)}:$ By Theorem \ref{Theorem IX}, the set of flat clopens of $\Min(A)$ forms a basis for the induced flat topology on $\Min(A)$. If $U\subseteq\Min(A)$ is a flat clopen, then by Theorem \ref{Theorem II1}, there exists an idempotent $e\in A$ such that $\gamma^{-1}(U)=V(e)$ where $\gamma:\Spec(A)\rightarrow\Min(A)$ is the retraction map. It follows that $U=V(e)\cap\Min(A)$. \\
$\mathbf{(vii)}\Rightarrow\mathbf{(i)}:$ If $\mathfrak{p}$ and $\mathfrak{q}$
are distinct minimal primes of $A$, then  $U=\Spec(A)\setminus\{\mathfrak{q}\}$ is an open neighborhood of $\mathfrak{p}$, because the closed points of $\Spec(A)$ with respect to the flat topology are precisely the minimal prime ideals of $A$. Thus by the hypothesis, there exists an idempotent $e\in A$ such that $\mathfrak{p}\in V(e)\cap\Min(A)\subseteq U\cap\Min(A)$. It follows that $e\in\mathfrak{p}$ and $1-e\in\mathfrak{q}$. \\
$\mathbf{(vii)}\Rightarrow\mathbf{(viii)}:$ If $I$ is a pure ideal of $A$, then $V(I)$ is a flat closed subset of $\Spec(A)$. Thus by the hypothesis, there exists a set $\{e_{i}\}$ of idempotents of $A$ such that $U\cap\Min(A)=\bigcup\limits_{i}\big(V(e_{i})\cap\Min(A)\big)$ where $U=\Spec(A)\setminus V(I)$. It follows that $U=\bigcup\limits_{i}V(e_{i})$, because $U$ is stable under the specialization. Therefore $V(I)=\bigcap\limits_{i}D(e_{i})=V(J)$ where the ideal $J$ is generated by the $1-e_{i}$. It follows that $I=J$, see \cite[Theorem 3.2]{Abolfazl 2020}. \\
$\mathbf{(viii)}\Rightarrow\mathbf{(iii)}:$ If $\mathfrak{p}$ is a minimal prime of $A$, then by Theorem \ref{Theorem VI}, it is a pure ideal. Thus by the hypotheses, $\mathfrak{p}$ is a regular ideal of $A$. \\
$\mathbf{(iv)}\Rightarrow\mathbf{(ix)}:$ If $M$ is a max-regular ideal of $A$, then $V(M)$ is a connected component of $\Spec(A)$. Thus by the hypothesis, there exists a minimal prime $\mathfrak{p}$ of $A$ such that $V(M)=V(\mathfrak{p})$. It follows that $M\subseteq\sqrt{M}=\mathfrak{p}$. But $\mathfrak{p}$ is a regular ideal of $A$, see (iii). This yields that $M=\mathfrak{p}$. By a similar argument, it is shown that every minimal prime ideal of $A$ is a max-regular ideal of $A$. \\
$\mathbf{(ix)}\Rightarrow\mathbf{(iv)}:$ For any ring $A$, the connected components of $\Spec(A)$ are precisely of the form $V(M)$ where $M$ is a max-regular ideal of $A$. \\
$\mathbf{(viii)}\Rightarrow\mathbf{(x)}:$ By Theorem \ref{Theorem VI}, $\Ann(f)$ is a pure ideal. Thus by the hypothesis, it is a regular ideal. \\
$\mathbf{(x)}\Rightarrow\mathbf{(xi)}:$ By Lemma \ref{Lemma nice 2}, there exists an idempotent $e\in\Ann(g)$ such that $f=fe$. It follows that $g=g(1-e)$. \\
$\mathbf{(xi)}\Rightarrow\mathbf{(i)}:$ Let $\mathfrak{p}$ and $\mathfrak{q}$ be distinct minimal prime ideals of $A$. By Lemma \ref{Lemma III}, there exist $f\in A\setminus\mathfrak{p}$ and $g\in A\setminus\mathfrak{q}$ such that $fg=0$. Then by the hypothesis, there is an idempotent $e\in A$ such that $f=fe$ and $g=g(1-e)$. This yields that $1-e\in\mathfrak{p}$ and $e\in\mathfrak{q}$. $\Box$ \\

The following result provides an elementwise characterization for reduced purified rings. \\

\begin{corollary}\label{Corollary Baer 9819} Let $A$ be a ring. Then $A$ is a reduced purified ring if and only if $fg=0$ for some $f,g\in A$, then there exists an idempotent $e\in A$ such that $f=fe$ and $g=g(1-e)$. \\
\end{corollary}

{\bf Proof.} For the implication ``$\Rightarrow$" see Theorem \ref{Theorem X} (xi). Conversely, if $f\in A$ is nilpotent then we may choose the least natural number $n\geqslant1$ such that $f^{n}=0$. If $n>1$ then there exists an idempotent $e\in A$ such that $f=fe$ and $f^{n-1}=f^{n-1}(1-e)$, but $f^{n-1}=f^{n-2}f(1-e)=0$ which is a contradiction. So $f=0$. Hence, $A$ is reduced. Thus by Theorem \ref{Theorem X} (xi), $A$ is also a purified ring. $\Box$ \\

\begin{corollary}\label{Corollary VIII} The direct product of a family of rings is a reduced purified ring if and only if each factor is a reduced purified ring. \\
\end{corollary}

{\bf Proof.} It follows from Corollary \ref{Corollary Baer 9819}. $\Box$ \\

There are reduced purified rings which are not p.p. rings, (recall that a ring $A$ is called a p.p. ring if each principal ideal of $A$ is a projective $A-$module, or equivalently, $\Ann(f)$ is generated by an idempotent element for all $f\in A$). For example, the ring $C(\beta\mathbb{N}\setminus\mathbb{N})$ is a reduced purified ring which is not a p.p. ring. There are also reduced mp-rings which are not purified rings. For instance, the ring $C(\beta\mathbb{R}^{+}\setminus\mathbb{R}^{+})$ is a reduced mp-ring which is not a purified ring. For the details see \cite{Al Ezeh 3}. Note that by Theorem \ref{Theorem X} (x), reduced purified rings and almost p.p. rings (in the sense of \cite{Al Ezeh 3}) are the same. Here, by $C(X)$ we mean the ring of real-valued continuous functions on $X$, and $\beta X$ denotes the Stone-\v{C}ech compactification of the topological space $X$. Hence, the following inclusions are strict: \\
p.p. rings $\subset$ reduced purified rings $\subset$ p.f. rings = reduced mp-rings. \\

\textbf{Acknowledgements.} The authors would like to give sincere thanks to Professors Pierre Deligne and Fran\c{c}ois Couchot for very valuable correspondences during the writing of the present paper. We would also like to give heartfelt thanks to the referees for very careful reading of the paper and for their very valuable comments and suggestions which improved the paper. \\

\end{document}